# A Hybrid Monte Carlo, Discontinuous Galerkin method for linear kinetic transport equations

Johannes Krotz[a], Cory D. Hauck[b], Ryan G. McClarren[c]

[a] *Department of Mathematics, University of Tennessee Knoxville, Knoxville, 37996, TN, USA*
[b] *Computational Mathematics Group, Computer Sciences and Mathematics Division, Oak Ridge National Laboratory, Oak Ridge, 37831, TN, USA*
[c] *Department of Aerospace and Mechanical Engineering, University of Notre Dame, Notre Dame, 46545, IN, USA*

**Abstract**

We present a hybrid method for time-dependent particle transport problems that combines Monte Carlo (MC) estimation with deterministic solutions based on discrete ordinates. For spatial discretizations, the MC algorithm computes a piecewise constant solution and the discrete ordinates uses bilinear discontinuous finite elements. From the hybridization of the problem, the resulting problem solved by Monte Carlo is scattering free, resulting in a simple, efficient solution procedure. Between time steps, we use a projection approach to "relabel" collided particles as uncollided particles. From a series of standard 2-D Cartesian test problems we observe that our hybrid method has improved accuracy and reduction in computational complexity of approximately an order of magnitude relative to standard discrete ordinates solutions.

## 1. Introduction

Numerical methods for kinetic transport equations are commonly divided into two classes: deterministic and Monte Carlo. Each of these approaches has strengths and weaknesses that complement the other.

Deterministic methods [23] directly discretize phase space (physical space, direction of flight, particle energy) as well as time, in the time-dependent setting. For this large seven-dimensional space (three for physical space, two for direction of flight, one for energy, and one for time), it is difficult to construct high resolution solutions for general problems. Indeed, the number of operations and the memory footprint required for deterministic solvers can be a challenge, even for leadership-class computers.

Monte Carlo methods [33], on the other hand, use sampling techniques to simulate particle transport processes. In its most basic form, the Monte-Carlo procedure is a computational analog of the actual physical processes being simulated: particles are sampled from sources and boundary conditions, then tracked as they stream through the domain, and along the way undergo scattering or absorption interactions with the material medium. As the particle traverses the physical domain, it contributes to integrated quantities of interest such as particle density or net fluence through a surface. For linear problems, the central limit theorem implies that the Monte Carlo solution is exact in the limit of an infinite number of samples [33]. Unlike deterministic methods, Monte Carlo methods are easy to extend to complicated 3-D geometries and can handle physical processes (such as particle interactions with the background material) in a continuous manner. Nevertheless, the uncertainty in Monte Carlo methods, as expressed in the standard deviation of an estimate, scales like $N^{-1/2}$, where $N$ is the number of sample particles. Additionally, Monte Carlo methods are not well-suited for obtaining uniform spatial estimates due to the difficulty of getting sufficient samples in every region of the physical domain. Moreover, for nonlinear problems such as thermal radiative transfer, the Monte Carlo approach loses some of its attractive properties. For example the discretization of material temperature, to which the particles are coupled, means that an exact solution is not obtained in the limit of infinite samples [42, 12, 25]. Nevertheless, producing efficient, accurate Monte Carlo calculations is an active area of research [32, 34].

Hybrid methods have been developed to harmonize the benefits of Monte Carlo and deterministic methods while minimizing their respective drawbacks. For steady-state nuclear reactor problems, methods such as COMET [26, 43] use local Monte Carlo calculations to estimate properties of solutions in macroscopic regions of the problem and then use a deterministic procedure to couple these regions together. Other work has considered weight windows and other biasing techniques [36, 6, 7, 37, 27] wherein deterministic solutions are used to modify the flight of particles in Monte Carlo calculations so that computational effort is spent more efficiently. High-order low-order (HOLO) schemes [28, 41, 29] have been developed in which Monte Carlo is used to solve compute a closure term for a low-order, moment-based deterministic calculation.

This work presents a deterministic-Monte Carlo hybrid method for time-dependent problems based on the physics of particle transport. Previous work [20, 9, 10, 21, 39] has exploited the fact describing a particle from its emission at the beginning of a time step to its first collision benefits from a different numerical treatment than a particle that emerges from a scattering interaction with the background medium. Because the scattering process relaxes particles towards a weakly anisotropic angular distribution, one can combine methods that are appropriate for particle streaming for the uncollided particles during a time step with methods that are suitable for weakly anisotropic angular distributions. In previous work, deterministic methods with a large number of angular degrees of freedom were used for the uncollided particles while low-resolution deterministic methods were used for the collided particles. It was also found that different levels of energy resolution can also be employed in this type of hybrid approach [40].

Despite the benefits of deterministic hybrid methods, solutions still require a large number of degrees of freedom for problems with large streaming paths. A natural strategy to address this challenge, which for steady-state problems was first proposed in [3], is to use Monte Carlo for the uncollided particles. Indeed, in many respects this is the ideal situation for a Monte Carlo approach. During a time step, particles are tracked through the computational domain and a non-analog estimator of the solution known as implicit capture is employed, thereby avoiding the need to consider collisions at all. Thus the calculation of the contribution to the solution from uncollided particles is essentially a ray tracing algorithm, which has many efficient implementations on modern computing hardware [2].

The hybrid method considered here uses Monte Carlo to compute the contribution to the solution from uncollided particles and an efficient deterministic calculation for the collided particles. A key advancement for extending the original steady-state formulation to time-dependent problems is a remapping step that resamples particles from the deterministic collided solution back into the uncollided component. This procedure is critical since, otherwise, the number of uncollided particles will decay exponentially and the hybrid solution will relax to a low-resolution, deterministic approximation of collided solution [20]. We find that the hybrid approach leads to more accurate solutions obtained using lower computational complexity than pure deterministic calculations.

The remainder of the paper is organized as follows. In Section 2, we introduce the hybrid method in the context of a single-group transport equation that is independent of particle energy. We also summarize the numerical methods used for the uncollided and collided components of the hybrid. In Section 3, we present numerical results for several standard test problems in a reduced two-dimensional geometry in physical space. In Section 4, we summarize findings and present directions for future work. A short appendix describes the Monte Carlo implementation of a boundary for one of the test problems.

## 2. Basics of the Hybrid Method

*2.1. Transport equation*

Let $X \in \mathbb{R}^3$ be a spatial domain with Lipschitz boundary and let $\mathbb{S}^2$ be the unit sphere in $\mathbb{R}^3$. Let $\Psi = \Psi(\boldsymbol{x}, \boldsymbol{\Omega}, t)$ be the angular flux depending on the position $\boldsymbol{x} = (x, y, z) \in X$, the direction of flight $\boldsymbol{\Omega} \in \mathbb{S}^2$ and time $t > 0$. We assume that $\Psi$ is governed by the linear transport equation

$$\frac{1}{c}\partial_t \Psi + \Omega \cdot \nabla_{\boldsymbol{x}} \Psi + \sigma_{\text{t}} \Psi = \frac{\sigma_{\text{s}}}{4\pi} \langle \Psi \rangle + Q, \qquad \boldsymbol{x} \in X, \quad \boldsymbol{\Omega} \in \mathbb{S}^2, \quad t > 0, \tag{1}$$

where $\sigma_{\text{t}} = \sigma_{\text{t}}(\boldsymbol{x})$, $\sigma_{\text{s}} = \sigma_{\text{s}}(\boldsymbol{x})$ and $\sigma_{\text{a}} = \sigma_{\text{t}} - \sigma_{\text{s}}$ are the total, scattering, and absorption cross-sections of the material, respectively; $Q = Q(\boldsymbol{x}, \boldsymbol{\Omega}, t)$ is a known particle source; and angle brackets denote integration over the unit sphere:

$$\langle \Psi \rangle = \int_{\mathbb{S}^2} \Psi \, d\boldsymbol{\Omega}. \tag{2}$$

The constant $c > 0$ is the particle speed; we assume that $c = 1$ for the remainder of this work. The transport equation (1) is equipped with initial conditions

$$\Psi(\boldsymbol{x}, \boldsymbol{\Omega}, 0) = \Psi_0(\boldsymbol{x}, \boldsymbol{\Omega}), \quad \boldsymbol{x} \in X, \quad \boldsymbol{\Omega} \in \mathbb{S}^2, \tag{3}$$

and boundary condition

$$\Psi(\boldsymbol{x}, \boldsymbol{\Omega}, t) = G(\boldsymbol{x}, \boldsymbol{\Omega}, t) \qquad \boldsymbol{x} \in \partial X, \quad \boldsymbol{\Omega} \cdot \boldsymbol{n}(\boldsymbol{x}) < 0, \tag{4}$$

where $\Psi_0$ and $G$ are known and $\boldsymbol{n}(\boldsymbol{x})$ is the unit outward normal at $\boldsymbol{x} \in \partial X$.

### 2.2. The hybrid method

The hybrid method is based on a first collision source splitting [3]. Let $\Psi = \Psi_u + \Psi_c$, where the *uncollided flux* $\Psi_u$ and the *collided flux* $\Psi_c$ satisfy the following system of equations

$$\partial_t \Psi_u + \boldsymbol{\Omega} \cdot \nabla_{\boldsymbol{x}} \Psi_u + \sigma_t \Psi_u = Q, \tag{5a}$$

$$\partial_t \Psi_c + \boldsymbol{\Omega} \cdot \nabla_{\boldsymbol{x}} \Psi_c + \sigma_t \Psi_c = \frac{\sigma_s}{4\pi} \left( \langle \Psi_u \rangle + \langle \Psi_c \rangle \right). \tag{5b}$$

Due to the linearity of (1), the splitting in (5) is exact; that is, if $\Psi_u$ and $\Psi_c$ solve (5a) and (5b), respectively, then $\Psi_u + \Psi_c$ solves (1). In practice, however, (5a) and (5b) are solved at different resolutions or even with different methods. Typically (5a) is solved with a method that has high resolution in angle, and because (5a) has no coupling in angle, it is easier to solve than (1) and also easy to solve in parallel. Meanwhile (5b) inherits the angular coupling in (1), but typically requires less angular resolution.

Since (5a) has no scattering source, particle densities will be transferred into the collided flux at an exponential rate, thus making the accuracy at which (5b) is solved the driving factor in the overall accuracy. This effect can be mitigated by abusing the autonomous nature of the equations and periodically relabeling the collided flux as uncollided at every time step.

To describe the implementation of the hybrid in more detail, let $\mathcal{T}$ be an operator such that

$$u[t] = \mathcal{T}(t, t', s, v, b, \lambda_t, \lambda_s) \tag{6}$$

where $u[t](\boldsymbol{x}, \boldsymbol{\Omega}) := u(\boldsymbol{x}, \boldsymbol{\Omega}, t)$, satisfies

$$\begin{cases} \partial_t u + \boldsymbol{\Omega} \cdot \nabla_{\boldsymbol{x}} u + \lambda_t u = \dfrac{\lambda_s}{4\pi} \langle u \rangle + s, & \boldsymbol{x} \in X, \quad \boldsymbol{\Omega} \in \mathbb{S}^2, \quad t > t', \tag{7a} \\ u(\boldsymbol{x}, \boldsymbol{\Omega}, t') = v(\boldsymbol{x}, \boldsymbol{\Omega}) & \boldsymbol{x} \in X, \quad \boldsymbol{\Omega} \in \mathbb{S}^2, \tag{7b} \\ u(\boldsymbol{x}, \boldsymbol{\Omega}, t) = b(\boldsymbol{x}, \boldsymbol{\Omega}, t) & \boldsymbol{x} \in \partial X, \quad \boldsymbol{\Omega} \cdot \boldsymbol{n}(\boldsymbol{x}) < 0, \quad t > t'. \tag{7c} \end{cases}$$

with source $s = s(\boldsymbol{x}, \boldsymbol{\Omega}, t)$. Using the operator $\mathcal{T}$, we can write

$$\Psi[t_{n+1}] = \mathcal{T}(t_{n+1}, t_n, Q, \Psi[t_n], G, \sigma_t, \sigma_s), \tag{8}$$

$$\Psi_u[t_{n+1}] = \mathcal{T}(t_{n+1}, t_n, Q_u, \Psi[t_n], G, \sigma_t, 0), \qquad Q_u = Q \tag{9}$$

$$\Psi_c[t_{n+1}] = \mathcal{T}(t_{n+1}, t_n, Q_c, 0, 0, \sigma_t, \sigma_s), \qquad Q_c = \frac{\sigma_s}{4\pi} \langle \Psi_u \rangle. \tag{10}$$

We simulate the system (5) using a Monte Carlo method for the uncollided equation (5a) and a deterministic discretization of the collided equation (5b). Let

$$\mathcal{T}_{\text{MC}}(t, t', s, v_{\text{MC}}, b, \lambda_t, \lambda_s; N_p) \tag{11}$$

be the Monte Carlo approximation to (7) given a particle representation $v_{\text{MC}}$ of $v$ and using $N_p$ pseudo-particles to represent the distribution of particles in phase space introduced by the source $s$ over the internal $(t', t)$. Similarly,

$$\mathcal{T}_{\text{SN}}(t, t', s, v, b, \lambda_t, \lambda_s; N, N_{\boldsymbol{x}}) \tag{12}$$

denote the $S_N$-DG approximation of (7) using a level $N$ set of ordinates, $N_{\boldsymbol{x}}$ spatial cells per dimension with $\mathbb{Q}_1$ elements, and a backward Euler time discretization to get from $t'$ to $t$. (The Monte Carlo method and SN-DG method are presented in more detail below.) Then given $N_p$, $N$, and $N_{\boldsymbol{x}}$, and a Monte Carlo approximation $\psi^n$ of $\Psi(t_n)$, let

$$\psi_{\mathrm{u}}^{n+1} = \mathcal{T}_{\mathrm{MC}}(t_{n+1}, t_n, Q_{\mathrm{u}}, \psi^n, G, \sigma_{\mathrm{t}}, 0; N_p), \qquad Q_{\mathrm{u}} = Q \tag{13}$$

$$\psi_{\mathrm{c}}^{n+1} = \mathcal{T}_{\mathrm{SN}}(t_{n+1}, t_n, Q_{\mathrm{c}}, 0, 0, \sigma_{\mathrm{t}}, \sigma_{\mathrm{s}}; N, N_{\boldsymbol{x}}), \qquad Q_{\mathrm{c}} = \frac{\sigma_{\mathrm{s}}}{4\pi} \langle \Psi_{\mathrm{u}} \rangle_{\mathrm{MC}} \tag{14}$$

$$\mathcal{R}\psi_{\mathrm{c}}^{n+1} = \mathcal{T}_{\mathrm{MC}}(t_{n+1}, t_n, Q_{\mathrm{r}}, 0, 0, \sigma_{\mathrm{t}}, 0; N_p), \qquad Q_{\mathrm{r}} = Q_{\mathrm{c}} + \frac{\sigma_{\mathrm{s}}}{4\pi} \langle \Psi_{\mathrm{c}} \rangle_{\mathrm{SN}} \tag{15}$$

$$\psi^{n+1} = \psi_{\mathrm{u}}^{n+1} + \mathcal{R}\psi_{\mathrm{c}}^{n+1} \tag{16}$$

where $\mathcal{R}$ is the remapping operator and $\langle \cdot \rangle_{\mathrm{MC}}$ and $\langle \cdot \rangle_{\mathrm{SN}}$ denote approximation of the angular integral over $\mathbb{S}^2$ with the respective method.

*2.3. Discrete ordinate-discontinuous Galerkin*

The discrete ordinates ($S_N$) method [5] approximates (7) by replacing the angular integral $\langle u \rangle$ by a discrete quadrature and then solving the resulting equation for the angles in the quadrature. This procedure yields a system of equations that depend only on space and time and can be further discretized by a variety of methods. Let

$$\{\boldsymbol{\Omega}_q\}_{q=1}^{N_{\boldsymbol{\Omega}}} \quad \text{and} \quad \{\omega_q\}_{q=1}^{N_{\boldsymbol{\Omega}}} \tag{17}$$

be the angles and associated weights of the $S_N$ quadrature, where $N_{\boldsymbol{\Omega}} = N_{\boldsymbol{\Omega}}(N)$ depends on the specific type of quadrature set being used. After discretizing in angle and applying an implicit Euler time discretization, the following semi-discrete system is obtained for each $q \in \{1, ..., N_{\boldsymbol{\Omega}}\}$ and $n \in \{0, 1, 2, \dots\}$,

$$\begin{cases} \dfrac{1}{\Delta t}\left(u_q^{n+1} - u_q^n\right) + \boldsymbol{\Omega}_q \cdot \nabla_{\boldsymbol{x}} u_q^{n+1} + \lambda_{\mathrm{t}} u_q^{n+1} = \dfrac{\lambda_{\mathrm{s}}}{4\pi} \displaystyle\sum_{r=1}^{N_{\boldsymbol{\Omega}}} \omega_r u_r^{n+1} + s_q^{n+1}, & \boldsymbol{x} \in X, \quad (18\mathrm{a}) \\ u_q^{n+1}(\boldsymbol{x}) = b_q^{n+1}(\boldsymbol{x}), & \boldsymbol{x} \in \partial X_q^-, \quad (18\mathrm{b}) \end{cases}$$

where $\partial X_q^- = \{\boldsymbol{x} \in X : \boldsymbol{\Omega}_q \cdot \boldsymbol{n}(\boldsymbol{x}) < 0\}$, $b_q^n(\boldsymbol{x}) = b(\boldsymbol{x}, \boldsymbol{\Omega}_q, t_n)$, $s_q^n(\boldsymbol{x}) = s(\boldsymbol{x}, \boldsymbol{\Omega}_q, t_n)$, and $u_q^n(\boldsymbol{x}) \approx u(\boldsymbol{x}, \boldsymbol{\Omega}_q, t_n)$ is the approximation on the temporal and angular grid. After reassigning

$$u_q \leftarrow u_q^{n+1}, \qquad s_q \leftarrow s_q^{n+1} + u_q^n, \qquad \lambda_{\mathrm{t}} \leftarrow \lambda_{\mathrm{t}} + \frac{1}{\Delta t}, \quad \text{and} \quad b_q \leftarrow b_q^{n+1}, \tag{19}$$

the discretization in (18a) can be written in the equivalent steady-state form

$$\begin{cases} \boldsymbol{\Omega}_q \cdot \nabla_{\boldsymbol{x}} u_q + \lambda_{\mathrm{t}} u_q = \lambda_{\mathrm{s}} \bar{u} + s_q, & \boldsymbol{x} \in X \quad (20\mathrm{a}) \\ u_q(\boldsymbol{x}) = b_q(\boldsymbol{x}), & \boldsymbol{x} \in \partial X_q^- \quad (20\mathrm{b}) \end{cases}$$

where $\mathbf{u} = (u_1, \dots, u_{N_{\boldsymbol{\Omega}}})^\top$, $\bar{u} := \frac{1}{4\pi} \sum_r \omega_r u_r$.

We discretize (20a) in physical space with a discontinuous Galerkin method and upwind numerical fluxes. The method by now is fairly standard (see for example [9, 19]) and is often used because of its accuracy in scattering-dominated regimes relative to upwind finite-difference and finite-volume methods [1, 22, 31, 18]. Because the DG method is well-known, we summarize it only briefly for the case of a two-dimensional Cartesian mesh with rectangular cells, which is sufficient for all of the numerical tests in Section 3. Let $X$ be divided into open sets $C_{i,j}$ that are squares with side lengths $h$ centered at $(x_i, y_j)$, and let $V_h = \{v \in L^2(X) : v|_{C_{ij}} \in \mathbb{Q}_1\}$. The goal will be to find an approximation of the weak solution of equation (20a); that is, find $\mathbf{u}^h = (u_1^h, \dots, u_{N_{\boldsymbol{\Omega}}}^h)^\top \in [V_h]^{N_{\boldsymbol{\Omega}}} := V_h \underbrace{\times \cdots \times}_{N_{\boldsymbol{\Omega}} \text{ times}} V_h$ such that

$$\mathcal{A}_q^{(i,j)}(u_q^h, v_q^h) + \mathcal{P}_q^{(i,j)}(u_q^h, v_q^h) = \mathcal{M}_q^{(i,j)}(u_q^h, v_q^h) + \mathcal{R}^{(i,j)}(\mathbf{u}^h, v_q^h) + \mathcal{S}_q^{(i,j)}(v_q^h) + \mathcal{B}_q^{(i,j)}(v_q^h) \tag{21}$$

for all $i,j \in \{1,\ldots,N_{\boldsymbol{x}}\}$, $q \in \{1,\ldots,N\}$, and $\mathbf{v}^h \in V_h^{N_{\boldsymbol{\Omega}}}$. Here

$$\mathcal{A}_q^{(i,j)}(u_q^h, v_q^h) = -\int_{C_{i,j}} (\boldsymbol{\Omega}_q \cdot \nabla_{\boldsymbol{x}} v_q^h) u_q^h d\boldsymbol{x} + \lambda_{\mathrm{t}} \int_{C_{i,j}} v_q^h u_q^h d\boldsymbol{x}, \tag{22}$$

$$\mathcal{P}_q^{(i,j)}(u_q^h, v_q^h) = \int_{(\partial C_{i,j})_q^+} (\boldsymbol{\Omega}_q \cdot \boldsymbol{n}) v_q^{h,-} v_q^{h,-} ds(\boldsymbol{x}), \qquad \mathcal{M}_q^{(i,j)}(u_q^h, v_q^h) = \int_{(\partial C_{i,j})_q^-} (\boldsymbol{\Omega}_q \cdot \boldsymbol{n}) v_q^{h,-} u_q^{h,+} ds(\boldsymbol{x}), \tag{23}$$

$$\mathcal{R}^{(i,j)}(\mathbf{u}^h, v_q^h) = \lambda_{\mathrm{s}} \int_{C_{i,j}} \bar{\mathbf{u}}^h v_q^h d\boldsymbol{x}, \quad \mathcal{S}_q^{(i,j)}(v_q^h) = \int_{C_{i,j}} s_q v_q^h d\boldsymbol{x}, \quad \mathcal{B}_q^{(i,j)}(v_q^h) = \int_{C_{i,j} \cap \partial X_q^-} b_q v_q^h ds(\boldsymbol{x}), \tag{24}$$

$$v_q^{h,\pm}(\boldsymbol{x}) = \lim_{\vartheta \to 0^+} v_q^h(\boldsymbol{x} \pm \vartheta \boldsymbol{n}), \quad \text{and} \quad (\partial C_{i,j})_q^{\pm} = \{\boldsymbol{x} \in \partial C_{i,j} : \pm \boldsymbol{\Omega}_q \cdot \boldsymbol{n}(\boldsymbol{x}) > 0\} \tag{25}$$

We construct $\mathbf{u}^h = (u_1^h, \ldots, u_N^h)^\top$ as follows: For each $i,j \in \{1,\ldots,N_{\boldsymbol{x}}\}$ and $q \in \{1,\ldots,N_{\boldsymbol{\Omega}}\}$, let

$$u_q^h(\boldsymbol{x}) = \sum_{|\mathbf{k}|_\infty \leq 1} \alpha_{q,\mathbf{k}}^{(i,j)} \phi_{\mathbf{k}}^{(i,j)}(x,y), \qquad \boldsymbol{x} \in C_{i,j}, \tag{26}$$

where $\mathbf{k} = (k_1, k_2)$,

$$\phi_{\mathbf{k}}^{(i,j)} = P_{k_1}\left(\frac{x - x_i}{h/2}\right) P_{k_2}\left(\frac{y - y_j}{h/2}\right), \tag{27}$$

and $P_k$ is the usual Legendre polynomial of degree $k$ on $\xi \in [-1, 1]$ with normalization $\int_{-1}^{1} P_{k_1}(\xi) P_{k_2}(\xi) d\xi = \frac{2}{2k_1+1}\delta_{k_1,k_2}$. Using this representation for $\mathbf{u}^h$, we derive the following linear system for the coefficients $\alpha_{q,\mathbf{k}}^{(i,j)}$ from equation (21): For each $\mathbf{l}$ such that $|\mathbf{l}|_\infty \leq 1$,

$$\sum_{|\mathbf{k}|_\infty \leq 1} \left(\mathbf{A}_{q,\mathbf{l},\mathbf{k}}^{(i,j)} + \mathbf{P}_{q,\mathbf{l},\mathbf{k}}^{(i,j)}\right) \alpha_{q,\mathbf{k}}^{(i,j)} = \sum_{|\mathbf{k}|_\infty \leq 1} \left(\mathbf{M}_{q,\mathbf{l},\mathbf{k}}^{(i,j)} \alpha_{q,\mathbf{k}}^{(i^*,j^*q)} + \mathbf{R}_{\mathbf{l},\mathbf{k}}^{(i,j)} \bar{\boldsymbol{\alpha}}_{\mathbf{k}}^{(i,j)}\right) + \mathbf{S}_{q,\mathbf{l}}^{(i,j)} + \mathbf{B}_{q,\mathbf{l}}^{(i,j)}, \tag{28}$$

where $\bar{\boldsymbol{\alpha}}_{\mathbf{k}}^{(i,j)} = \sum_{q=1}^{N_{\boldsymbol{\Omega}}} w_q \alpha_{q,\mathbf{k}}^{(i,j)}$,

$$\mathbf{A}_{q,\mathbf{l},\mathbf{k}}^{(i,j)} = \mathcal{A}_q^{(i,j)}(\phi_{\mathbf{k}}^{(i,j)}, \phi_{\mathbf{l}}^{(i,j)}), \qquad \mathbf{P}_{q,\mathbf{l},\mathbf{k}}^{(i,j)} = \mathcal{P}_q^{(i,j)}(\phi_{\mathbf{k}}^{(i,j)}, \phi_{\mathbf{l}}^{(i,j)}) \qquad \mathbf{M}_{q,\mathbf{l},\mathbf{k}}^{(i,j)} = \mathcal{M}_q^{(i,j)}(\phi_{\mathbf{k}}^{(i^*,j^*)}, \phi_{\mathbf{l}}^{(i,j)}) \tag{29}$$

$$\mathbf{B}_{q,\mathbf{l}}^{(i,j)} = \mathcal{B}_q^{(i,j)}(\phi_{\mathbf{l}}^{(i,j)}), \qquad \mathbf{S}_{q,\mathbf{l}}^{(i,j)} = \mathcal{S}_q^{(i,j)}(\phi_{\mathbf{l}}^{(i,j)}), \qquad \mathbf{R}_{\mathbf{l},\mathbf{k}}^{(i,j)} = \mathcal{R}^{(i,j)}(\phi_{\mathbf{k}}^{(i,j)}, \phi_{\mathbf{l}}^{(i,j)}) \tag{30}$$

and, given the components $\boldsymbol{n} = (n_x, n_y)$ of the outward normal, the indices $i^*, j^*$ are given by

$$i^*(\boldsymbol{x}) = i + n_x(\boldsymbol{x}) \quad \text{and} \quad j^*(\boldsymbol{x}) = j + n_y(\boldsymbol{x}). \tag{31}$$

To improve readability, we rewrite equation (28) as a matrix equation with respect to the indices $\mathbf{k}$ and $\mathbf{l}$:

$$\left(\mathbf{A}_q^{(i,j)} + \mathbf{P}_q^{(i,j)}\right) \boldsymbol{\alpha}_q^{(i,j)} = \mathbf{R}^{(i,j)} \bar{\boldsymbol{\alpha}}^{(i,j)} + \mathbf{M}_q^{(i,j)} \boldsymbol{\alpha}_q^{(i^*,j^*)} + \mathbf{B}_q^{(i,j)} + \mathbf{S}_q^{(i,j)}. \tag{32}$$

The organization of the operators in (32) reflects a standard solution strategy combining source iteration and sweeping. In source iteration, $\bar{\boldsymbol{\alpha}}^{(i,j)}$ is lagged; that is, given an iteration index $\ell$:

$$\boldsymbol{\alpha}_q^{(i,j,\ell+1)} = \left(\mathbf{A}_q^{(i,j)} + \mathbf{P}_q^{(i,j)}\right)^{-1} \left(\mathbf{R}^{(i,j)} \bar{\boldsymbol{\alpha}}^{(i,j,\ell)} + \mathbf{M}_q^{(i,j)} \boldsymbol{\alpha}_q^{(i^*,j^*,\ell+1)} + \mathbf{B}_q^{(i,j)} + \mathbf{S}_q^{(i,j)}\right), \tag{33a}$$

$$\bar{\boldsymbol{\alpha}}^{(i,j,\ell+1)} = \sum_{q=1}^{N} w_q \boldsymbol{\alpha}_q^{(i,j,\ell+1)}. \tag{33b}$$

Sweeping refers to process solving of (33) cell-by-cell: for each $q$, cells can be ordered such that $\boldsymbol{\alpha}_q^{(i^*,j^*,\ell+1)}$ is known, prior to solving for $\boldsymbol{\alpha}_q^{(i,j,\ell+1)}$. The details of this procedure are given in Algorithm 1.

*2.4. Monte Carlo*

In this section, we describe the Monte Carlo method used to compute the solution to (7) for the pure absorption problem when $\lambda_t = \lambda_a$ (i.e. no scattering):

$$\begin{cases} \partial_t u + \mathbf{\Omega} \cdot \nabla_{\boldsymbol{x}} u + \lambda_t u = s, & \boldsymbol{x} \in X, \quad \mathbf{\Omega} \in \mathbb{S}^2, \quad t > 0, & \text{(34a)} \\ u(\boldsymbol{x}, \mathbf{\Omega}, 0) = v(\boldsymbol{x}, \mathbf{\Omega}) & \boldsymbol{x} \in X, \quad \mathbf{\Omega} \in \mathbb{S}^2, & \text{(34b)} \\ u(\boldsymbol{x}, \mathbf{\Omega}, t) = b(\boldsymbol{x}, \mathbf{\Omega}, t) & \boldsymbol{x} \in \partial X, \quad \mathbf{\Omega} \cdot \boldsymbol{n}(\boldsymbol{x}) < 0, \quad t > 0. & \text{(34c)} \end{cases}$$

The Monte Carlo method approximates the phase space distribution $u$ using a finite set of pseudo-particles:

$$u(\boldsymbol{x}, \mathbf{\Omega}, t) \approx \sum_{\pi \in \Pi^t} w_\pi(t) \delta(\boldsymbol{x} - \boldsymbol{x}_\pi(t)) \delta(\mathbf{\Omega} - \mathbf{\Omega}_\pi), \tag{35}$$

where $\Pi^t$ is a set of pseudo-particles $\pi$ with position $x_\pi(t)$, weight $w_\pi(t)$, and direction of flight $\mathbf{\Omega}_\pi$. $\Pi^t$ will be defined more carefully below. A benefit of the hybrid is that scattering processes, which can slow down the method significantly [13], do not need to be modeled in (34).

The Monte Carlo implementation of (34) can be derived via a Green's function formulation for (34). Let $G(\boldsymbol{x}, \boldsymbol{y}, \mathbf{\Omega}, t, t_0)$ solve

$$\partial_t G + \mathbf{\Omega} \cdot \nabla_x G + \lambda_t G = \delta(\boldsymbol{x} - \boldsymbol{y}) \delta(t - t_0), \tag{36}$$

with zero initial data and boundary conditions. Then

$$u(\boldsymbol{x}, \mathbf{\Omega}, t) = \int_0^t \int_{\mathbb{R}^3} G(\boldsymbol{x}, \boldsymbol{y}, \mathbf{\Omega}, t, t_0) s(\boldsymbol{y}, \mathbf{\Omega}, t_0) d\boldsymbol{y} dt_0 \tag{37a}$$

$$+ \int_0^t \int_{\mathbb{R}^3} G(\boldsymbol{x}, \boldsymbol{y}, \mathbf{\Omega}, t, t_0) s_v(\boldsymbol{y}, \mathbf{\Omega}, t_0) d\boldsymbol{y} dt_0 \tag{37b}$$

$$+ \int_0^t \int_{\mathbb{R}^3} G(\boldsymbol{x}, \boldsymbol{y}, \mathbf{\Omega}, t, t_0) s_b(\boldsymbol{y}, \mathbf{\Omega}, t_0) d\boldsymbol{y} dt_0 \tag{37c}$$

solves (34), where $s$, $s_v$, and $s_b$ are identically zero outside of the closure of $X$. The terms $s_v$ and $s_b$ are provisional source terms designed such that (37b) solves (34) when $s = 0$ and $b = 0$, while (37c) solves (34) when $s = 0$ and $v = 0$. The former is solved by setting $s_v(\boldsymbol{x}, \mathbf{\Omega}, t) = v(\boldsymbol{x}, \mathbf{\Omega}) \delta(t)$. However, determining $s_b$ can be slightly more involved, and an example that is used for numerical experiments in Section 3.3 is provided in the Appendix. Once $s, s_v$, and $s_b$ are known, all three terms can be treated identically. For simplicity, we restrict our attention to (37a) below.

For any $t > t_0$ and any fixed $\mathbf{\Omega}$, let $\boldsymbol{X}(t) = \boldsymbol{x}_0 + (t - t_0)\mathbf{\Omega}$. Then $g(t) = G(\boldsymbol{X}(t), \boldsymbol{y}, \mathbf{\Omega}, t, t_0)$ solves

$$\frac{dg(t)}{dt} = -\lambda(\boldsymbol{X}(t))g(t) + \delta(\boldsymbol{X}(t) - \boldsymbol{y})\delta(t - t_0) \tag{38}$$

or, equivalently,

$$g(t) = e^{-\int_{t_0}^t \lambda_t(\boldsymbol{X}(\xi))d\xi} \delta(\boldsymbol{x}_0 - \boldsymbol{y}). \tag{39}$$

Setting $\boldsymbol{x} = \boldsymbol{X}(t)$ in (39) gives

$$G(\boldsymbol{x}, \boldsymbol{y}, \mathbf{\Omega}, t, t_0) = e^{-\int_{t_0}^t \lambda_t(\boldsymbol{x} - (t-\xi)\mathbf{\Omega})d\xi} \delta(\boldsymbol{x} - (t - t_0)\mathbf{\Omega} - \boldsymbol{y}). \tag{40}$$

Plugging (40) back into (37a) yields, after some manipulation,

$$\begin{aligned}
u(\boldsymbol{x},\boldsymbol{\Omega},t) &= \int_0^t \int_{\mathbb{R}^3} e^{-\int_{t_0}^t \lambda_{\mathrm{t}}(\boldsymbol{x}-(t-\xi)\boldsymbol{\Omega})d\xi} \delta(\boldsymbol{x}-(t-t_0)\boldsymbol{\Omega}-\boldsymbol{y})s(\boldsymbol{y},\boldsymbol{\Omega},t_0)d\boldsymbol{y}dt_0 \\
&= \int_0^t e^{-\int_{t_0}^t \lambda_{\mathrm{t}}(\boldsymbol{x}-(t-\xi)\boldsymbol{\Omega})d\xi} s(\boldsymbol{x}-(t-t_0)\boldsymbol{\Omega},\boldsymbol{\Omega},t_0)dt_0 \\
&= \int_0^t e^{-\int_{t-\tau}^t \lambda_{\mathrm{t}}(\boldsymbol{x}-(t-\xi)\boldsymbol{\Omega})d\xi} s(\boldsymbol{x}-\tau\boldsymbol{\Omega},\boldsymbol{\Omega},t-\tau)d\tau \\
&= \int_0^t e^{-\int_0^\tau \lambda_{\mathrm{t}}(\boldsymbol{x}-\xi\boldsymbol{\Omega})d\xi} s(\boldsymbol{x}-\tau\boldsymbol{\Omega},\boldsymbol{\Omega},t-\tau)d\tau.
\end{aligned} \qquad (41)$$

This representation of $u(\boldsymbol{x},\boldsymbol{\Omega},t)$ can be interpreted as the density of particles that have reached the location $\boldsymbol{x}$ at time $t$ while moving in the direction $\boldsymbol{\Omega}$. These particles are emitted by the source $s$ at time $t-\tau$ and location $\boldsymbol{x}-t\boldsymbol{\Omega}$, and they carry a weight that decays exponentially due to absorption.

The Monte Carlo approach can be understood as an approximation of $u$ based on sampling of pseudo-particles from the source $s$ in (41). Let

$$\tilde{s}(\boldsymbol{x},\boldsymbol{\Omega},t) = \sum_{\pi \in \Pi} w_\pi \delta(\boldsymbol{x}-\boldsymbol{x}_\pi)\delta(\boldsymbol{\Omega}-\boldsymbol{\Omega}_\pi)\delta(t-t_\pi) \qquad (42)$$

where $\pi \in \Pi$ are pseudo-particles with weight $w_\pi > 0$, position $\boldsymbol{x}_\pi \in X$, and direction of flight $\boldsymbol{\Omega}_\pi \in \mathbb{S}^2$ at $t_\pi > 0$, such that $\tilde{s}(\boldsymbol{x},\boldsymbol{\Omega},t)$ approximates $s(\boldsymbol{x},\boldsymbol{\Omega},t)$ for all $t > 0$. Then
for any $C \subset X$, any $B \subset \mathbb{S}^2$, and any time interval $(t_n,t_{n+1})$, the representation of $u$ in (41), along with the approximation $\tilde{s}$ gives

$$\begin{aligned}
&\int_{t_n}^{t_{n+1}} \int_C \int_B u(\boldsymbol{x},\boldsymbol{\Omega},t)d\boldsymbol{\Omega} d\boldsymbol{x}dt \\
&\approx \int_{t_n}^{t_{n+1}} \int_C \int_B \int_0^t e^{-\int_0^\tau \lambda_{\mathrm{t}}(\boldsymbol{x}-\xi\boldsymbol{\Omega})d\xi} \tilde{s}(\boldsymbol{x}-\tau\boldsymbol{\Omega},\boldsymbol{\Omega},t-\tau)d\tau d\boldsymbol{\Omega} d\boldsymbol{x}dt \\
&= \int_{t_n}^{t_{n+1}} \int_C \int_B \int_0^t e^{-\int_0^\tau \lambda_{\mathrm{t}}(\boldsymbol{x}-\xi\boldsymbol{\Omega})d\xi} \sum_{\pi \in \Pi} w_\pi \delta(\boldsymbol{x}-\tau\boldsymbol{\Omega}-\boldsymbol{x}_\pi)\delta(\boldsymbol{\Omega}-\boldsymbol{\Omega}_\pi)\delta(t-\tau-t_\pi)d\tau d\boldsymbol{\Omega} d\boldsymbol{x}dt \\
&= \sum_{\pi \in \Pi} \int_{t_n}^{t_{n+1}} w_\pi e^{-\int_0^{t-t_\pi} \lambda_{\mathrm{t}}(\boldsymbol{x}_\pi+\xi\boldsymbol{\Omega}_\pi)d\xi} \mathbb{1}_C(\boldsymbol{x}_\pi+(t-t_\pi)\boldsymbol{\Omega}_\pi)\mathbb{1}_{[0,t]}(t_\pi)\mathbb{1}_B(\Omega_\pi)dt \\
&= \sum_{\pi \in \Pi} \int_{t_n}^{t_{n+1}} w_\pi(t)\mathbb{1}_C(\boldsymbol{x}_\pi+(t-t_\pi)\boldsymbol{\Omega}_\pi)\mathbb{1}_{[0,t]}(t_\pi)\mathbb{1}_B(\Omega_\pi)dt
\end{aligned} \qquad (43)$$

where $w_\pi(t) = w_\pi e^{-\int_0^{t-t_\pi} \lambda_{\mathrm{t}}(\boldsymbol{x}_\pi+\xi\boldsymbol{\Omega}_\pi)d\xi}$ and $w_\pi(t_\pi) = w_\pi$.

Note that with the identification $\boldsymbol{x}_\pi(t) = \boldsymbol{x}_\pi + (t-t_\pi)\boldsymbol{\Omega}_\pi$ we can also identify $\Pi^t$ in (35) as $\Pi^t = \{\pi \in \Pi : t_\pi < t\}$.

We will denote (43) as $\mathcal{T}_{\mathrm{MC}}(t_{n+1},t_n,0,\lambda_{\mathrm{s}},\lambda_{\mathrm{s}},s,0;N_p)$ as the Monte Carlo solution for the case of the zero initial data and zero boundary data. A general Monte-Carlo solution can be obtained as

$$\begin{aligned}
\mathcal{T}_{\mathrm{MC}}(t_{n+1},t_n,s,v,b,\lambda_{\mathrm{s}},\lambda_{\mathrm{s}};N_p) &= \mathcal{T}_{\mathrm{MC}}(t_{n+1},t_n,s+s_b,v,0,\lambda_{\mathrm{s}},\lambda_{\mathrm{s}};N_p) \\
&= \mathcal{T}_{\mathrm{MC}}(t_{n+1},t_n,s+s_v+s_b,0,0,\lambda_{\mathrm{s}},\lambda_{\mathrm{s}};N_p)
\end{aligned}$$

Numerically it is useful to realize that

$$\Pi^{n+1} := \Pi^{t_{n+1}} = \{\pi(t_n+\Delta t) : \pi \in \Pi^{t_n}\} \cup \{\pi \in \Pi : t_\pi \in (t_n,t_{n+1})\} \qquad (44)$$

where in an abuse of notation $\pi(t_n+\Delta t)$ is a new particle with position $x_\pi(t+\Delta t)$, weight $w_\pi(t+\Delta t)$ and direction of flight $\boldsymbol{\Omega}_\pi$ for some $\pi \in \Pi^{t_n} = \Pi^n$. Thus particles $\Pi^{n+1}$ can be obtained by updating the weight and position of particles in $\Pi^n$ and sampling new particles with birth times in $(t_n, t_{n+1})$. From (43) we also obtain

$$\langle u \rangle_{\text{MC}}(\boldsymbol{x}, t) = \sum_{\pi \in \Pi} \int_{t_n}^{t_{n+1}} w_\pi(t) \mathbb{1}_C(\boldsymbol{x}_\pi + (t - t_\pi)\boldsymbol{\Omega}_\pi) \mathbb{1}_{[0,t]}(t_\pi) dt \qquad (45)$$

With this formulation, particle weights $w_\pi(t)$ decrease exponentially at a rate proportional to the absorption, given by $\lambda_t$. This approach is known as the implicit capture method. It avoids the need to sample absorption times explicitly and, at the same time, reduces statistical noise and simplifies the implementation [14, Page 168][24, Chapter 22].

The Monte Carlo simulation of (34) from $t_n$ to $t_{n+1}$ is based on (43) and proceeds according to the following steps:

1. Let $\mathcal{P}$ be a partition of $X$ into disjoint cells $C$. For each $C \in \mathcal{P}$, calculate the total weight $W^C$ of new particles generated in $C$ by the source during the interval $(t_n, t_{n+1})$:

$$W^C = \frac{1}{\Delta t} \frac{1}{4\pi} \int_{\mathbb{S}^2} \int_{t_n}^{t_{n+1}} \int_C s(\boldsymbol{x}, \Omega, t) \, d\Omega d\boldsymbol{x} dt, \qquad (46)$$

where $\Delta t = t^{n+1} - t^n$, and let $W = W^C \sum_{C \in \mathcal{P}}$. Let $N_p$ be the input for the total number of new particles to sample during the interval $(t_n, t_{n+1})$. Then for each $C \in \mathcal{P}$, sample

$$N_p^C = \text{floor}\left[\frac{W^C}{W}\right] \qquad (47)$$

   particles and assign them each a weight $w = W/N_p^C$. The total number of particles in the system at this time is $\tilde{N}_p^{n+1} = \tilde{N}_p^{n+1} + \sum_{C \in \mathcal{P}} N_p^C$. Each new particle $p$ is assigned a position $\boldsymbol{x}_\pi$ sampled uniformly from $C$ and a birth time $t^{n+1} - \tau_\pi$, where $\tau_\pi$ is sampled uniformly from $(0, \Delta t)$. Each particle $p$ is also assigned an angle $\boldsymbol{\Omega}_\pi$. For isotropic sources, $\boldsymbol{\Omega}_\pi$ is sampled uniformly from $\mathbb{S}^2$. For non-isotropic sources, (such as the boundary source for the holhraum problem in Section 3.3), the sampling distribution must be consistent with the angular dependence. The particles, including their space, angle, and time coordinates, are added to current particle list.

2. Move each particle $\pi$ in the current particle list from $\boldsymbol{x}_\pi$ to $x_\pi + \tau_\pi \boldsymbol{\Omega}_\pi$ and update its weight to $w_\pi \leftarrow w_\pi(t_{n+1})$. The number of particles in the system will have to be adjusted accordingly. Reset its remaining time to $\tau_\pi \leftarrow \Delta t$. For each cell $C \in \mathcal{P}$, update the sum in (43) and (45) according the time spent in $C$ during the interval $(t_n, t_{n+1})$.

3. To reduce computational effort and memory, at the end of each time step any particle $p$ with weight $w_\pi < w_{\text{kill}}$ will be dropped with a probability of $p_{\text{kill}} > 0$. Here $w_{\text{kill}} > 0$ is a user-defined parameter, called the 'killing weight', and $p_{\text{kill}} = (1 - w_\pi/w_{\text{kill}})$. To preserve the total mass in the system, any particle $p$ with weight $w < w_{\text{kill}}$ that survives this 'Russian roulette' [24, Chapter 22] will have its weight readjusted to $w_\pi/(1 - p_{\text{kill}})$.

*2.5. Pseudocode*

The algorithms that we use for our numerical results are detailed in Algorithms 1-3. The $S_N$ method is given in Algorithm 1. The hybrid method is given in Algorithm 3. It requires Algorithm 1 for the collided component and Algorithm 2 for the Monte Carlo update. A listing of the notation used these algorithms is provided in Table 1.

## 3. Numerical Results

In this section, we compare simulation results from the Monte Carlo-$S_N$ hybrid method to those from a standard, monolithic $S_N$ method. The goal is to demonstrate that the hybrid method provides a more efficient approach. The $S_N$ computations for the monolithic method and for the collided component of the hybrid rely on product quadrature sets on the sphere [38, 4].

| User defined parameters | |
|---|---|
| $N_{\boldsymbol{x}}$ | number of Cartesian cells along each dimension |
| $N$ | order of discrete ordinates |
| $N_p$ | number of new particles (up to rounding) generated from source |
| $\Delta t$ | time step |
| $\delta$ | tolerance of iteration |
| $\omega_q$ | Gauss-Legendre weights |
| **Material parameters** | |
| $\lambda_{\mathrm{t}}, \lambda_{\mathrm{a}}, \lambda_{\mathrm{s}}$ | total, absorption and scattering crosssection |
| **Additional notation** | |
| $\mathbf{A}_q^{(i,j)}, \mathbf{P}_q^{(i,j)}, \mathbf{R}^{(i,j)}, \mathbf{M}_q^{(i,j)}, \mathbf{B}_q^{(i,j)}, \mathbf{S}_q^{(i,j)}$ | matrices defined in Section 2.3 |
| $\mathcal{U}(X)$ | uniform distribution on set $X$. |

Table 1: Pseudocode parameters

---

**Algorithm 1** $S_N$-algorithm: Propagate solution from $t_k$ to $t_{k+1}$

**Input:** $u^n, \mathbf{s}_q$   ▷ coefficients of solution $u_q^k$ from previous step and source
**Require:** $\hat{\lambda}_t, \lambda_{\mathrm{a}}, \lambda_{\mathrm{s}},$   ▷ Material properties
**Require:** $\Delta t, \{C_{i,j}\}_{i,j=1}^{N_{\boldsymbol{x}}}, \{\boldsymbol{\Omega}_q, w_q\}_{q=1}^{N_{\boldsymbol{\Omega}}}$   ▷ Discretization parameters
**Require:** $\delta$   ▷ Convergence tolerance

1: $\boldsymbol{\alpha}_q^{(i,j)}(t_n)$ such that $u_q^n(\boldsymbol{x}) = \sum_{|\mathbf{k}|_\infty \leq 1} \alpha_{q,\mathbf{k}}^{(i,j)} \phi_{\mathbf{k}}^{(i,j)}(\boldsymbol{x})$ for $\boldsymbol{x} \in C_{i,j}$
2: $\mathbf{s}_q \leftarrow \mathbf{s}_q + \frac{1}{\Delta t} \boldsymbol{\alpha}_q(t_n)$   ▷ Initialize source
3: $\boldsymbol{\alpha}_q \leftarrow \boldsymbol{\alpha}_q(t_n)$   ▷ Initialize coefficients
4: $err = \delta + 1$
5: **while** $err > \delta$ **do**
6:     $\boldsymbol{\beta}_q \leftarrow \boldsymbol{\alpha}_q$   ▷ Store old coefficients
7:     $\bar{\boldsymbol{\alpha}} \leftarrow \sum_{q=1}^{N_{\boldsymbol{\Omega}}} w_q \boldsymbol{\alpha}_q$
8:     **for** $q \in \{1, ..., N_{\boldsymbol{\Omega}}\}$ **do**
9:         **for** $(i,j) \in \{1, ..., N_{\boldsymbol{x}}\}^2$ **do**   ▷ Sweep through cells in direction $\boldsymbol{\Omega}_q$
10:             $\boldsymbol{\alpha}_q^{(i,j)} \leftarrow \left( \mathbf{A}_q^{(i,j)} + \mathbf{P}_q^{(i,j)} \right)^{-1} \left( \mathbf{R}^{(i,j)} \bar{\beta}^{(i,j)} + \mathbf{M}_q^{(i,j)} \boldsymbol{\alpha}_q^{(i^*, j^*)} + \mathbf{B}_q^{(i,j)} + \mathbf{S}_q^{(i,j)} \right)$
  ▷ Update coefficients
11:         **end for**
12:     **end for**
13:     $err = \max_q \| \boldsymbol{\alpha}_q - \boldsymbol{\beta}_q \|_\infty$   ▷ Discrepancy between iterations
14: **end while**
15: **return** $u_q^{n+1}(\boldsymbol{x}) = \sum_{|\mathbf{k}|_\infty \leq 1} \alpha_{q,\mathbf{k}}^{(i,j)} \phi_{\mathbf{k}}^{(i,j)}(\boldsymbol{x})$, $\langle u_q^{n+1} \rangle_{\mathrm{SN}}(\boldsymbol{x}) = \sum_{|\mathbf{k}|_\infty \leq 1} \bar{\alpha}_{\mathbf{k}}^{(i,j)} \phi_{\mathbf{k}}^{(i,j)}(\boldsymbol{x})$ for $\boldsymbol{x} \in C_{i,j}$

### Algorithm 2 MC-algorithm: Propagate solution form $t_k$ to $t_{k+1}$

---

**Input:** $\Pi^{n+1}$ with $u^n(\boldsymbol{x}, \boldsymbol{\Omega}) = \sum_{\pi \in \Pi^t} w_\pi \delta(\boldsymbol{x} - \boldsymbol{x}_\pi)\delta(\boldsymbol{\Omega} - \boldsymbol{\Omega}_\pi)$

$s(\boldsymbol{x}, \boldsymbol{\Omega}, t)$  ▷ previous particle distribution and source

**Require:** $\lambda_t$

**Require:** $\Delta t, \{C_{i,j}\}_{i,j}^{N_x}, N_p$  ▷ Discretization parameters

**Require:** $w_{\text{kill}}$  ▷ killing weight

1: **for** $(i,j) \in \{1, \ldots, N_x\}^2$ **do**
2: $\quad W^{C_{i,j}} = \frac{1}{|C_{i,j}|\Delta t} \int_{t_n}^{t_{n+1}} \int_{C_{i,j}} \int_{\mathbb{S}^2} s(\boldsymbol{x}, \boldsymbol{\Omega}, t) d\boldsymbol{\Omega} d\boldsymbol{x} dt$
3: **end for**
4: $W \leftarrow \sum_{i,j} W^{C_{i,j}}$  ▷ Calculate total source
5: **for** $(i,j) \in \{1, \ldots, N_x\}^2$ **do**
6: $\quad N_p^{C_{i,j}} \leftarrow \left\lfloor \frac{W^{C_{i,j}} N_p}{W} \right\rfloor$  ▷ Number of particles on each cell
7: **end for**
8: $N_p \leftarrow \sum_{i,j} N_p^{C_{i,j}}$  ▷ number of new particles
9: $w \leftarrow \frac{W}{N_p}$
10: **for** $(i,j)$ with $N_p^{C_{i,j}} > 0$ **do**  ▷ sample new particles from source
11: $\quad$ **for** $k \in \{1, \ldots, N_p^{C_{i,j}}\}$ **do**
12: $\quad\quad$ Generate new particle $\pi$ with
13: $\quad\quad \boldsymbol{x}_\pi \sim \mathcal{U}(C_{i,j})$  ▷ Draw particle's position
14: $\quad\quad (\Omega_x, \Omega_y, \Omega_z) \sim \frac{1}{W^{C_{i,j}}|C_{i,j}|\Delta t} \int_{t_n}^{t_{n+1}} \int_{C_{i,j}} s(\boldsymbol{x}, \boldsymbol{\Omega}, t) d\boldsymbol{x} dt$  ▷ Draw particle's direction of flight
15: $\quad\quad \boldsymbol{\Omega}_\pi \leftarrow (\Omega_x, \Omega_y)$
16: $\quad\quad w_\pi \leftarrow w$  ▷ Assign particle weight
17: $\quad\quad \Pi^* \leftarrow \Pi^* \cup \{\pi\}$  ▷ Add new particle to existing
18: $\quad$ **end for**
19: **end for**
20: **for** $\pi \in \Pi^t \cup \Pi^*$ **do**  ▷ Move particles
21: $\quad$ **if** $\pi \in \Pi^t$ **then**
22: $\quad\quad \tau_\pi \leftarrow \Delta t$  ▷ remaining time for particles from prev. step
23: $\quad$ **else**
24: $\quad\quad$ randomly draw $\tau_\pi \sim \mathcal{U}([0, \Delta t])$  ▷ remaining time for particles generated this time step
25: $\quad$ **end if**
26: $\quad \boldsymbol{x}_\pi \leftarrow \boldsymbol{x}_\pi + \tau_\pi \boldsymbol{\Omega}_\pi$  ▷ Update particle's position
27: $\quad$ **for** $(i,j)$ with $C_{i,j} \cap \{\boldsymbol{x}\pi + t\boldsymbol{\Omega}_\pi : t \in [0, \tau_\pi]\} \neq \emptyset$ **do**  ▷ All cells intersected by particle's trajectory
28: $\quad\quad \Phi_{i,j} \leftarrow \Phi_{i,j} + w_\pi \int_0^{\tau_\pi} \exp\left(-\int_0^t \lambda_a(\boldsymbol{x}_p + t'\boldsymbol{\Omega}_p)dt'\right) \mathbb{1}_{C_{i,j}}(\boldsymbol{x}_\pi + t\boldsymbol{\Omega}_\pi)dt$  ▷ Update $\Phi$
29: $\quad$ **end for**
30: $\quad w_\pi \leftarrow w_\pi \exp\left(-\int_0^{\tau_\pi} \lambda_a(\boldsymbol{x}_\pi + t\boldsymbol{\Omega}_p)dt\right)$  ▷ Update particle weight
31: $\quad$ **if** $\boldsymbol{x}_\pi \in X$ **then**
32: $\quad\quad \Pi^{n+1} \leftarrow \Pi^{n+1} \cup \{\pi\}$  ▷ Remove particles that left domain
33: $\quad$ **end if**
34: **end for**
35: **for** $\pi \in \Pi^{n+1}$ with $w_\pi < w_{\text{kill}}$ **do**  ▷ Russian roulette
36: $\quad r \sim \mathcal{U}([0,1])$
37: $\quad$ **if** $r > \frac{w_\pi}{w_{\text{kill}}}$ **then**  ▷ Determine survival of particle
38: $\quad\quad \Pi^{n+1} \leftarrow \Pi^{n+1} \setminus \{\pi\}$
39: $\quad$ **else**
40: $\quad\quad w_\pi \leftarrow w_{\text{kill}}$  ▷ Update surviving particle's weight to approx. preserve total mass
41: $\quad$ **end if**
42: **end for**
43: **return** $\Pi^{n+1}$ with $u^{n+1}(\boldsymbol{x}, \boldsymbol{\Omega}) = \sum_{\pi \in \Pi^{n+1}} w_\pi \delta(\boldsymbol{x} - \boldsymbol{x}_\pi)\delta(\boldsymbol{\Omega} - \boldsymbol{\Omega}_\pi)$ and $\langle u^{n+1} \rangle_{\text{MC}}(\boldsymbol{x}) = \Phi(\boldsymbol{x})$

---

**Algorithm 3** $H_{\text{MC-S}_\text{N}}$-algorithm: Propagate solution from $t_k$ to $t_{k+1}$
---
    **Input:** $\Pi^n$ with $u^n(\boldsymbol{x}, \boldsymbol{\Omega}) = \sum_p w_p \delta(\boldsymbol{x} - \boldsymbol{x}_p) \delta(\boldsymbol{\Omega} - \boldsymbol{\Omega}_p)$

    $s(\boldsymbol{x}, \boldsymbol{\Omega}, t)$      ▷ previous particle distribution and source

**Require:** $\hat{\lambda}_\text{t}, \lambda_\text{t}, \lambda_\text{a}, \lambda_\text{s},$      ▷ Material properties

**Require:** $\Delta t, \{C_{i,j}\}_{i,j=1}^{N_x}, \{\boldsymbol{\Omega}_q, w_q\}_{q=1}^{N_\Omega}, N_p$      ▷ Discretization parameters

**Require:** $\delta, w_{\text{kill}}$      ▷ Convergence tolerance and killing weight

1:  $\Pi_u^{n+1}, \langle u_u^{n+1}\rangle_{\text{MC}} \leftarrow \text{MC}(\Pi^n, s)$      ▷ Monte Carlo for uncollided
2:  **for** $q \in \{1, ..., N_\Omega\}$ **do**
3:     $s_q \leftarrow \lambda_\text{s} \Phi_u$      ▷ turn $\langle u_u^{n+1}\rangle_{\text{MC}}$ into source for $S_n$
4:  **end for**
5:  $u_c^{n+1}, \langle u_c^{n+1}\rangle_{\text{SN}} \leftarrow S_n(0, s_q)$      ▷ $S_n$ for collided
6:  $s \leftarrow \lambda_\text{s}(\langle u_c^{n+1}\rangle_{\text{SN}} + \langle u_u^{n+1}\rangle_{\text{MC}})$      ▷ sources for Relabeling
7:  $\Pi_R, \langle u_R^{n+1}\rangle_{\text{MC}} \leftarrow \text{MC}(0, s)$      ▷ Monte Carlo as relabeling
8:  $\Pi^{n+1} \leftarrow \Pi_u^{n+1} \cup \Pi_R$
9:  $\langle u^{n+1}\rangle_{\text{MC}} \leftarrow \langle u_u^{n+1}\rangle_{\text{MC}} + \langle u_R^{n+1}\rangle_{\text{MC}}$
10: **return** $\Pi^{n+1}, \langle u^{n+1}\rangle_{\text{MC}}$

---

We consider three well known test problems: the line source problem [16], the lattice problem [8], and the linearized hohlraum problem [9, 8]. The specifications for each problem are provided in the following subsections. They are all formulated in a geometry for which $\partial_z \Psi = 0$. This means that they can be reduced to two dimensions in physical space and, by an abuse of notation, we write $\Psi(\boldsymbol{x}, \boldsymbol{\Omega}, t) = \Psi(x, y, \boldsymbol{\Omega}, t)$. A further consequence of the geometry is that product quadrature on the sphere can be reduced to just the upper hemisphere, in which case $N_\Omega = N^2$. The time step for each problem is tied to the grid resolution via the ratio $\text{CFL} = \Delta t / \Delta x$. For all calculations shown below, the iteration tolerance $\delta$ (see Algorithm 1) is set to $10^{-4}$.[1]

For each problem, we assess the accuracy of the numerical solution and the efficiency with which it is obtained. To quantify the accuracy we compare our results to a reference solution. For the line source problem the reference is the semi-analytic solution from [17]; see also [16]. For the other two test problems, the reference is a high-resolution hybrid solution based solely on $S_N$ discretization with a triangular-based quadrature referred to as $T_N$ [35] for both collided and uncollided component, combined with a DG discretization in space and integral deferred correction in time [11].

Accuracy for the MC-$S_N$ hybrid and the monolithic $S_N$ method is measured in terms of the relative $L^2$-difference in the scalar flux $\Phi = \langle \Psi \rangle$ at a given final time $t_{\text{final}}$. Given the numerical solution $\Phi_{\text{num}}$ and the reference $\Phi_{\text{ref}}$:

$$\Delta = \frac{||\Phi - \Phi_{\text{ref}}||_{L_2}}{||\Phi_{\text{ref}}||_{L_2}}, \tag{48}$$

where $L^2$-norm is approximated by $h^2 \sum_{C_{i,j}} \Phi_{ij}^2$ and $\Phi_{i,j}$ is the average on the cell $C_{i,j}$. Because our implementation of the hybrid method and the $S_N$ method are not run-time optimized, we use a complexity measure which counts the number of times a particle is moved or a *DG* unknown is updated in the course of a sweep. Let $N^c$ be the level of the quadruature for the collided component of the hybrid. Then the complexity of the monolithic method is

$$\mathbb{C}_{S_N} = \underset{\substack{\uparrow \\ \text{\# of Legendre} \\ \text{coefficients}}}{4} \times \underset{\substack{\uparrow \\ \text{\# of or-} \\ \text{dinates}}}{N_\Omega} \times \underset{\substack{\uparrow \\ \text{\# of} \\ \text{cells}}}{N_x^2} \times \underset{\substack{\uparrow \\ \text{\# of source} \\ \text{iterations}}}{N_i} \times \underset{\substack{\uparrow \\ \text{\# of time} \\ \text{steps}}}{\frac{T}{\Delta t}} \tag{49a}$$

---
[1] While not shown here, we also tested several runs using $\delta = 10^{-8}$, which leads to negligible improvements in accuracy when compared to the results shown below.

while the complexity of the hybrid is

$$\mathbb{C}_{\text{hybrid}} = \underbrace{(N_u}_{\substack{\text{avg. \# of} \\ \text{particles} \\ \text{moved for} \\ \text{uncollided}}} + \underbrace{N_R)}_{\substack{\text{avg. \#} \\ \text{particles} \\ \text{moved in} \\ \text{relabelling}}} \times \underbrace{\frac{T}{\Delta t}}_{\substack{\text{\# of time} \\ \text{steps}}} + \underbrace{\mathbb{C}_{S_{Nc}}}_{\substack{S_N \\ \text{complexity of} \\ \text{collided} \\ \text{equation}}} \qquad (49\text{b})$$

For convenience, we set $N_{\text{MC}}^{\text{tot}} = (N_u + N_R)\frac{T}{\Delta t}$ so that $\mathbb{C}_{\text{hybrid}} = N_{\text{MC}}^{\text{tot}} + \mathbb{C}_{S_{Nc}}$. $N_u$ is the sum of particles added to the system $N_p$ and the average number of particles still in flight from the previous time step $N_{prev}$, i.e. $N_u = N_p + N_{prev}$, while $N_R$ is the number of particles added in the relabeling, which we set to be $N_R = N_p$.

### 3.1. The Line Source problem

In the line source problem, an initial pulse of uniformly distributed particles is emitted from the line $\ell = \{(x, y, z) : x = y = 0\}$ into the surrounding domain $X = \mathbb{R}^3$, which contains a purely scattering material with $\sigma_s = \sigma_t = 1$. Because the geometry of the domain and initial condition are invariant in $z$, the spatial domain can be reduced to $\mathbb{R}^2$. In this two-dimensional setting, the initial condition can be represented by an isotropic delta function $\frac{1}{4\pi}\delta(x, y)$, but to reduce numerical artifacts, we use a mollified version of the initial condition:

$$\Psi_0(x, y, \mathbf{\Omega}, t) = \frac{1}{4\pi} \frac{1}{2\pi\varsigma} \exp\left(-\frac{x^2 + y^2}{2\varsigma}\right), \qquad \varsigma = 0.03. \qquad (50)$$

Meanwhile, the computational domain is restricted to the square $[-1.5, 1.5]^2$ and equipped with zero inflow boundary conditions.

We perform monolithic $S_N$ and MC-$S_N$ hybrid simulations at various spatial and angular resolutions. The spatial domain is subdivided into equal $N_{\boldsymbol{x}} \times N_{\boldsymbol{x}}$ square cells with $N_{\boldsymbol{x}} \in \{51, 101, 201\}$. For the $S_N$-runs we let $N \in \{4, 8, 16, 32\}$, resulting in $N_{\boldsymbol{\Omega}} = N^2$ ordinates on the northern hemisphere of $\mathbb{S}^2$. The collided part of the hybrid algorithm employs an $S_N$ method with $N \in \{4, 8\}$. In the hybrid method the number of particles was also changed between runs. The number of particles newly inserted into the system every time step is roughly $2 \times N_p$, where $N_p = 10^k$ and $k \in \{2, 3, 4, 5, 6\}$.[2] Due to rounding and particles being dropped via Russian roulette, the exact number of particles inserted into the system varies slightly. The killing weight is fixed at $w_{\text{kill}} = 10^{-15}$. The CFL is also fixed at 0.5 across all runs. The reference solution is the semi-analytic solution from [17]; see also [16].

Figure 1 depicts several approximations of the scalar flux $\Phi$ at $t_{\text{final}} = 1$, calculated using the $S_N$ method and the hybrid method. The solutions calculated using the $S_N$ method clearly show ray-effects that only get resolved after a significant increase in the angular resolution. No such effects are seen in the hybrid solutions. The hybrid solutions do contain some noise, as the particle count is relatively low, but they preserve the symmetry of the problem up to a reasonable error. Unlike the $S_N$-method, the hybrid is able to capture the wave front of unscattered particles travelling away from the center with speed 1.

Figure 2 shows $L_2$-errors of the numerical solutions in log scale; the same trends are apparent. While the hybrid solutions mainly suffer from noise due to the stochastic nature of the Monte Carlo method, the $S_N$ method has strong ray-effects and struggles to capture the analytic solution at the wave front.

A more systematic analysis of the numerical results is presented in Figure 3. This plot shows the relative $L_2$-error $\Delta$ of various runs versus their respective computational complexity $\mathbb{C}$. For the hybrid method, increasing the angular resolution $N$ in the collided component yields a marginal improvement at best. However, changes in the particle number $N_p$ for the uncollided component have a significant impact. For the $S_N$ method, on the other hand, increasing the angular resolution yields a significant improvement in the accuracy. For smaller values of $N$, increasing the spatial resolution may actually increase the error. This is especially apparent in Figure 3 for the $S_4$ and $S_8$ results. For a fixed angular resolution, additional spatial accuracy will begin to resolve the ray effect anomalies in the solution. Conversely, $S_N$ results with lower spatial resolution benefit from error cancellation due to the numerical diffusion smoothing ray effects. The hybrid may also have larger errors if the particles per cell is too low.

---

[2]The factor of 2 is because $N_p$ particles are used for the uncollided equations (13) and $N_p$ particles are used for the relabeling (15).

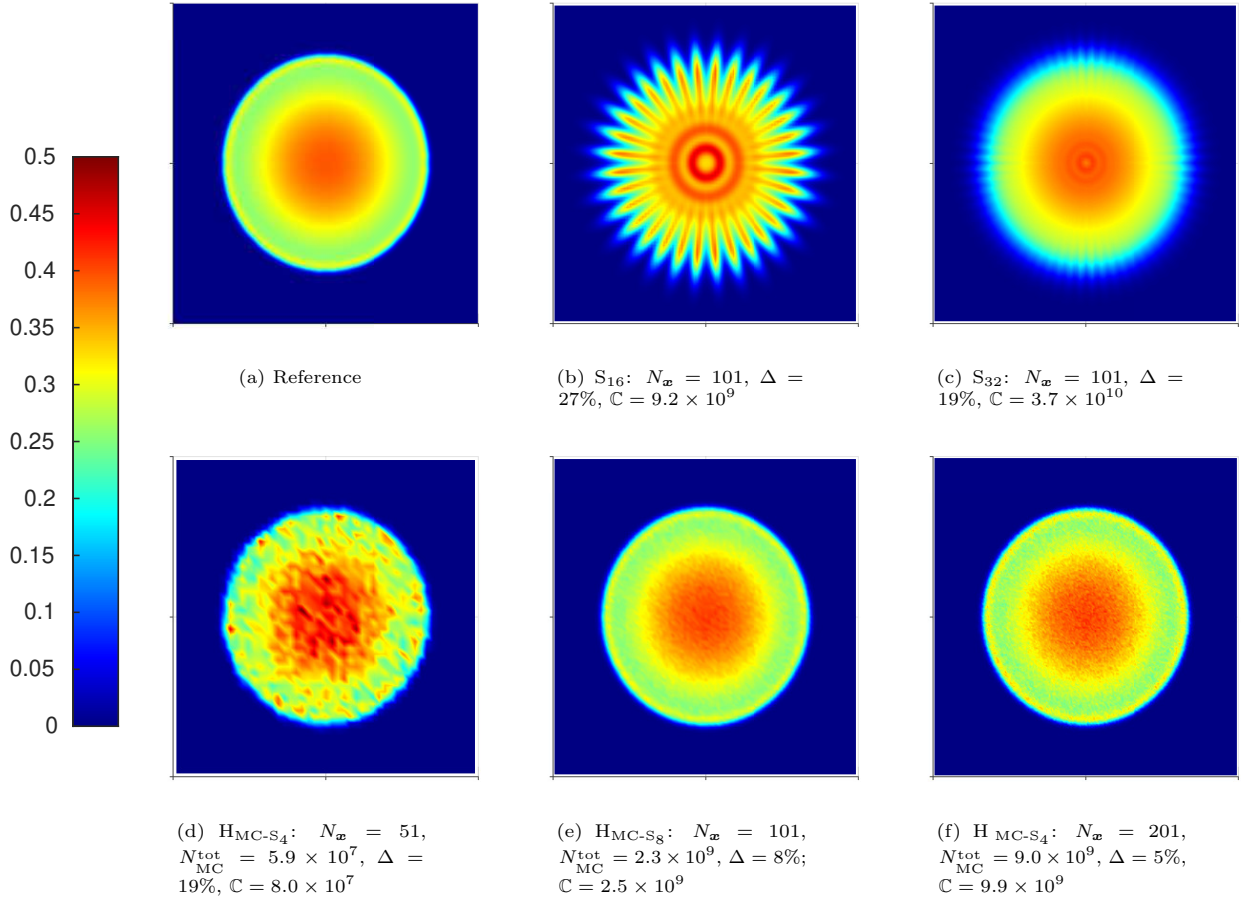

Figure 1: Numerical approximation of the scalar flux $\Phi$ for the line source problem at $t_{\text{final}} = 1$ with CFL 0.5. Each numerical solution is characterized by a relative $L^2$ difference $\Delta$ with respect to the reference, defined in (48), and a complexity $\mathbb{C}$, defined in (49a) for the monolithic $S_N$ method and (49b) for the hybrid.

Overall the hybrid method outperforms the monolithic $S_N$ method. For example, the hybrid error can match the most refined $S_N$ calculation ($N = 32$) with a complexity that is roughly 2-3 orders of magnitude smaller; compare Figures 2 (c) and (d). In fact the hybrid method can obtain an error with half the size with a complexity that is an order of magnitude less. In general hybrid runs tend to be 3-4 times more accurate than their $S_N$ counterparts of similar complexity.

*3.2. The Lattice problem*

In the lattice problem, a checkerboard of highly absorbing material is embedded in a scattering material with a central source. The layout of this problem along with its material parameters can be found in Figure 4. The computational domain is a $7 \times 7$ rectangle with zero inflow data at the boundaries. The center square (red) contains an isotropic particle source, while the blue squares are pure absorbers. The red and white squares are purely scattering with $\sigma_{\text{s}} = \sigma_{\text{t}} = 1$. The initial condition is identically zero everywhere in the domain.

We perform $S_N$ and hybrid runs with varying spatial and angular resolution. The spatial domain is subdivided into equal $N_{\boldsymbol{x}} \times N_{\boldsymbol{x}}$ square cells with $N_{\boldsymbol{x}} \in \{56, 112, 224\}$. For the $S_N$-runs we use $N \in \{4, 8, 16, 32\}$, resulting in $N_{\boldsymbol{\Omega}} = N^2$ ordinates on the northern hemisphere of $\mathbb{S}^2$. The collided part of the hybrid algorithm employs an $S_N$ method with $N \in \{4, 8\}$. In the hybrid method the number of particles is also changed between runs. The number of particles newly inserted into the system every time step is roughly $2 \times N_p$ where $N_p = 10^k$ for $k \in \{2, 3, 4, 5, 6\}$. The killing weight is fixed at $w_{\text{kill}} = 10^{-15}$. All runs are performed to a final time of $t_{\text{final}} = 3.2$ and the CFL is kept fixed at 25.6. The reference solution is a $S_{96}$-$S_{16}$ hybrid (meaning a $S_{96}$ for the uncollided component and $S_{16}$ for the collided component) using a $T_N$ quadrature in angle, a third-order DG method in space on a 448 by 448 grid, and a defect correction time integrator [10].

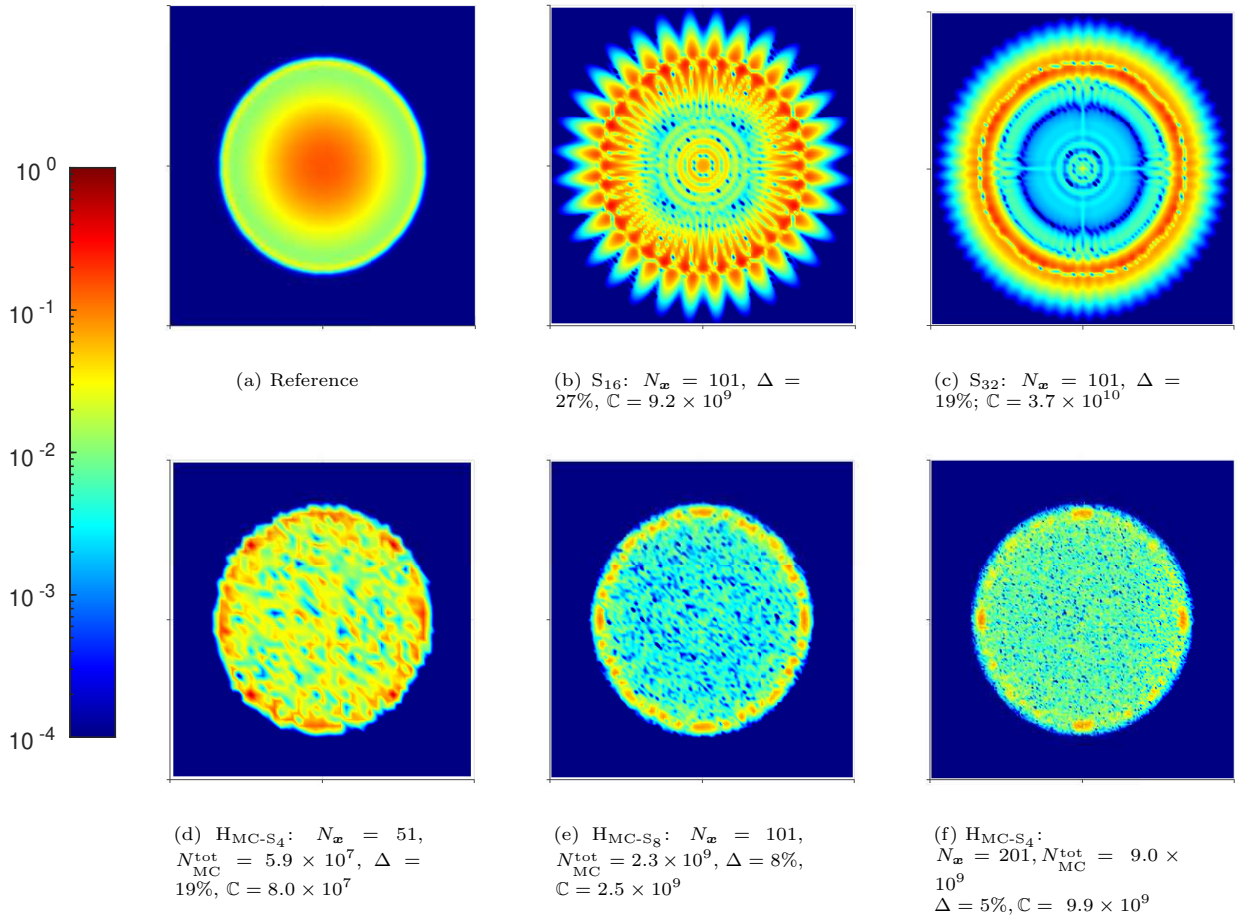

Figure 2: Absolute difference between the analytical solution and various numerical solutions to the line source problem at $t = 1$ with CFL 0.5. Each numerical solution is characterized by a relative $L^2$ difference $\Delta$ with respect to the reference, defined in (48), and a complexity $\mathbb{C}$, defined in (49a) for the monolithic $S_N$ method and (49b) for the hybrid.

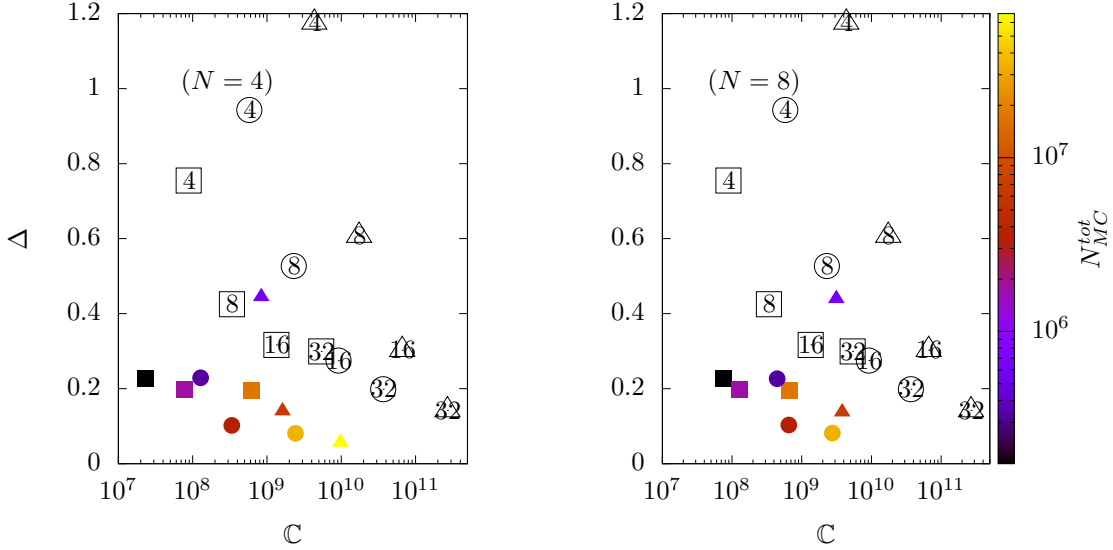

Figure 3: The relative $L^2$-difference $\Delta$ vs. complexity $\mathbb{C}$ for the scalar flux $\Phi$ in the line source problem, using the hybrid method (filled, colored markers) and the monolithic $S_N$ method (empty, green markers). Points that are down and to the left are more efficient. The hybrid solver was run for $N = 4$ (left) and $N = 8$ (right). The shape of data-points indicates different spatial resolution (square: $N_{\boldsymbol{x}} = 51$, circle: $N_{\boldsymbol{x}} = 101$, triangle: $N_{\boldsymbol{x}} = 201$). Coloring of the hybrid data points corresponds to the total number of particles per time step $N_{MC}^{tot}$ according to the colorbar. The $S_N$ method was run for $N = 4, 8, 16, 32$. Each DG-data point is assigned a numerical label according to its value of $N$. The formula for $\Delta$ is given in (48) which the complexity $\mathbb{C}$ is given by (49a) for the monolithic $S_N$ method and by (49b) for the hybrid.

Selected results for the scalar flux $\Phi$ are depicted in Figure 5, and the relative $L_2$-error for these same solutions is depicted in Figure 6. While the hybrid solutions are not completely free of ray-effects, these effects are much more pronounced in the $S_N$ runs. A more rigorous analysis of the performance of the two algorithm in dependence of their respective parameters can be seen in Figure 7. This plot shows the relative error $\Delta$ of various runs in dependence of their complexity $\mathbb{C}$. It is noteworthy that for this test problem the accuracy is mostly independent of the angular resolution, but depends significantly on the spatial resolution. Increasing the overall particle count in the hybrid method is most effective at higher spatial resolution; compare for example the vertical separation in colored triangles vs. colored circles vs. colored squares in Figure 7.

It turns out that increasing the number of particles in the hybrid algorithm does not necessarily increase the algorithm complexity. This is because with increased particle count the iterative solver for the collided components often needs fewer iterations. In cases where the complexity is dominated by these iterations, an increase in particles can even cause a decrease in computational complexity. Overall, Figure 7 shows that the hybrid algorithm produces results with comparable or slightly better accuracy than the standard $S_N$ solver, while being close to an order of magnitude of lower complexity.

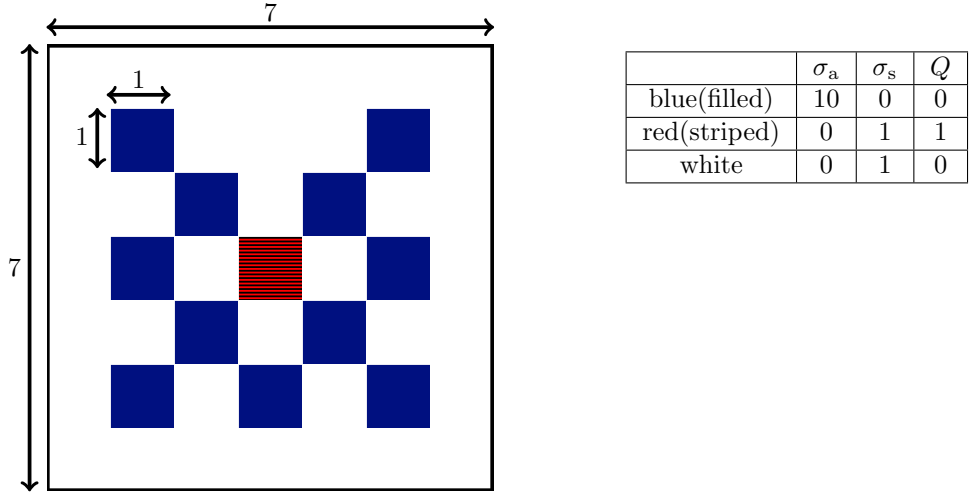

Figure 4: Geometric layout and table of material properties for the Lattice Problem.

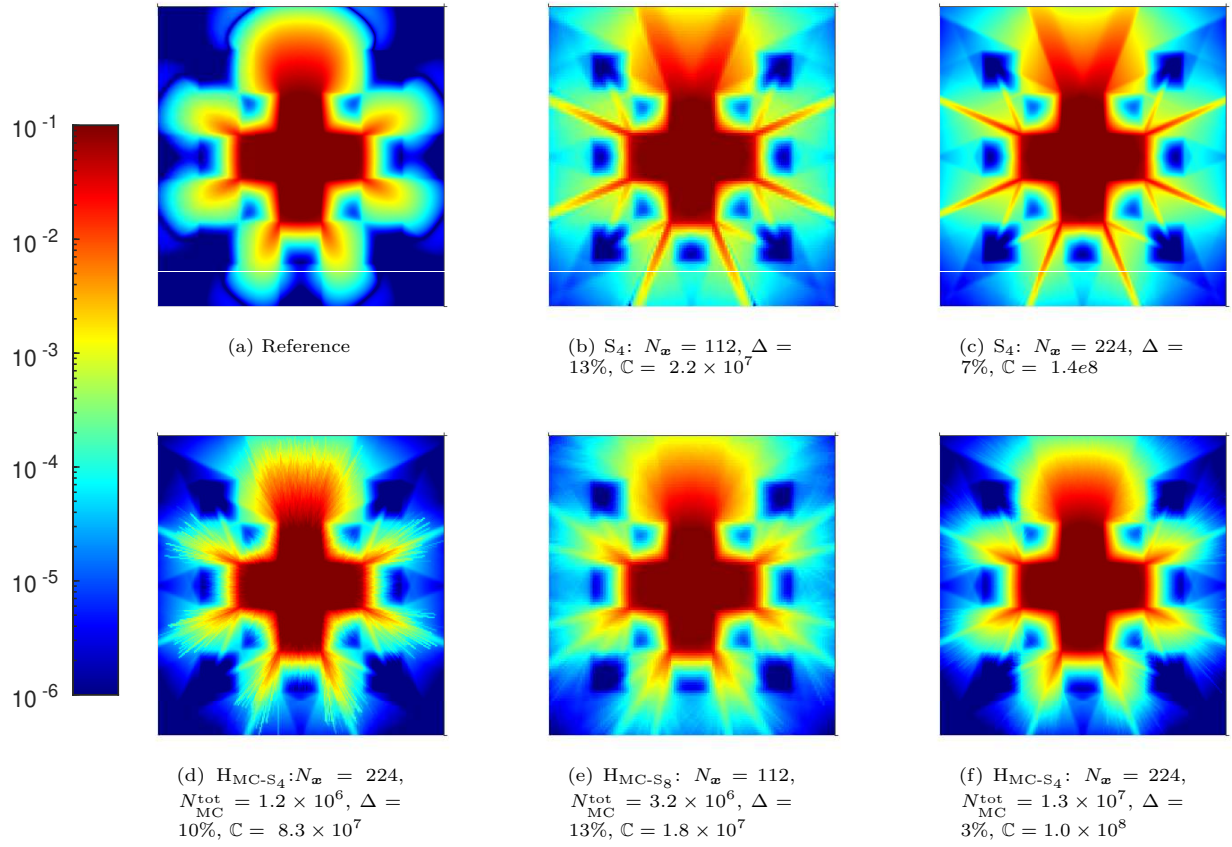

(a) Reference

(b) $S_4$: $N_x = 112$, $\Delta = 13\%$, $\mathbb{C} = 2.2 \times 10^7$

(c) $S_4$: $N_x = 224$, $\Delta = 7\%$, $\mathbb{C} = 1.4e8$

(d) $\text{H}_{\text{MC-}S_4}$: $N_x = 224$, $N_{\text{MC}}^{\text{tot}} = 1.2 \times 10^6$, $\Delta = 10\%$, $\mathbb{C} = 8.3 \times 10^7$

(e) $\text{H}_{\text{MC-}S_8}$: $N_x = 112$, $N_{\text{MC}}^{\text{tot}} = 3.2 \times 10^6$, $\Delta = 13\%$, $\mathbb{C} = 1.8 \times 10^7$

(f) $\text{H}_{\text{MC-}S_4}$: $N_x = 224$, $N_{\text{MC}}^{\text{tot}} = 1.3 \times 10^7$, $\Delta = 3\%$, $\mathbb{C} = 1.0 \times 10^8$

Figure 5: Numerical solutions to the lattice problem at $t = 3.2$ with CFL 25.6. Each numerical solution is characterized by a relative $L^2$ difference $\Delta$ with respect to the reference, defined in (48), and a complexity $\mathbb{C}$, defined in (49a) for the monolithic $S_N$ method and (49b) for the hybrid.

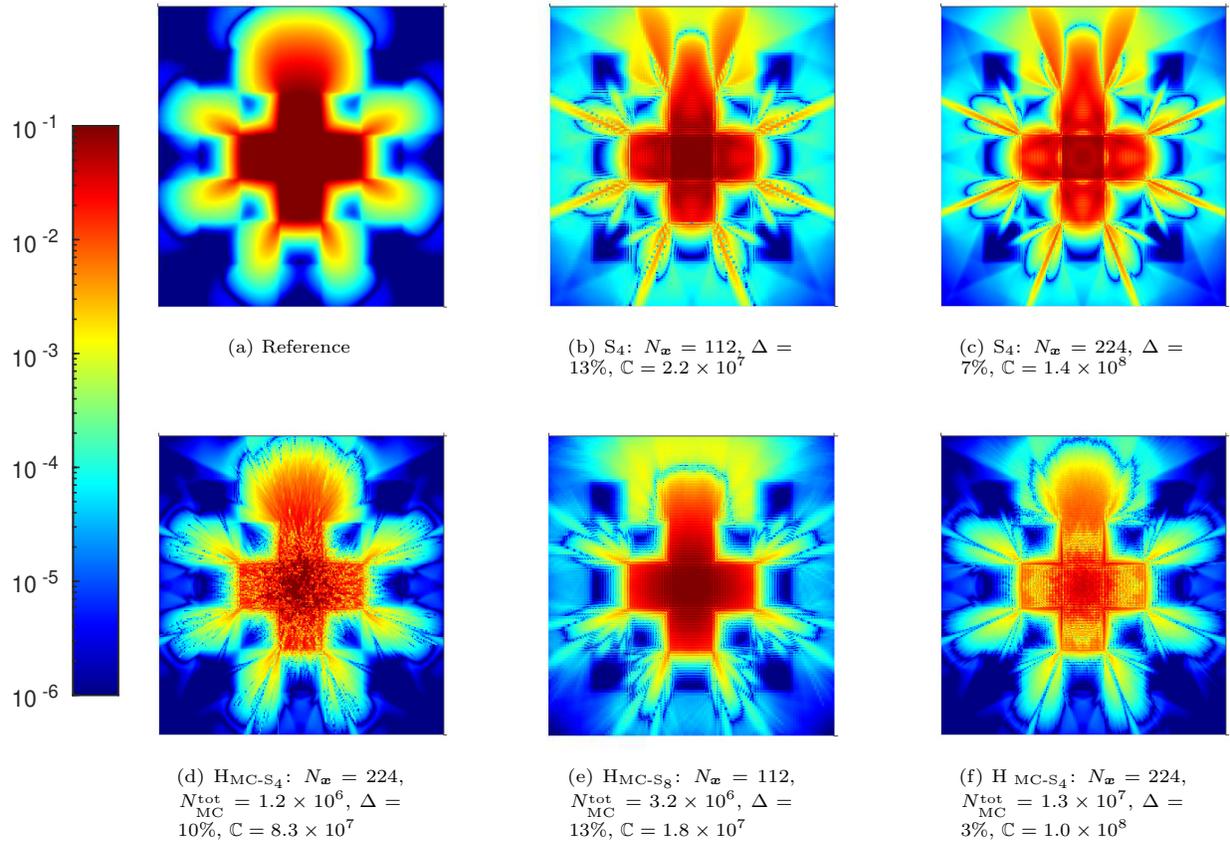

Figure 6: Lattice problem: Absolute difference between the reference solution and various numerical solutions at $t = 3.2$ with CFL 25.6. Each numerical solution is characterized by a relative $L^2$ difference $\Delta$ with respect to the reference, defined in (48), and a complexity $\mathbb{C}$, defined in (49a) for the monolithic $S_N$ method and (49b) for the hybrid.

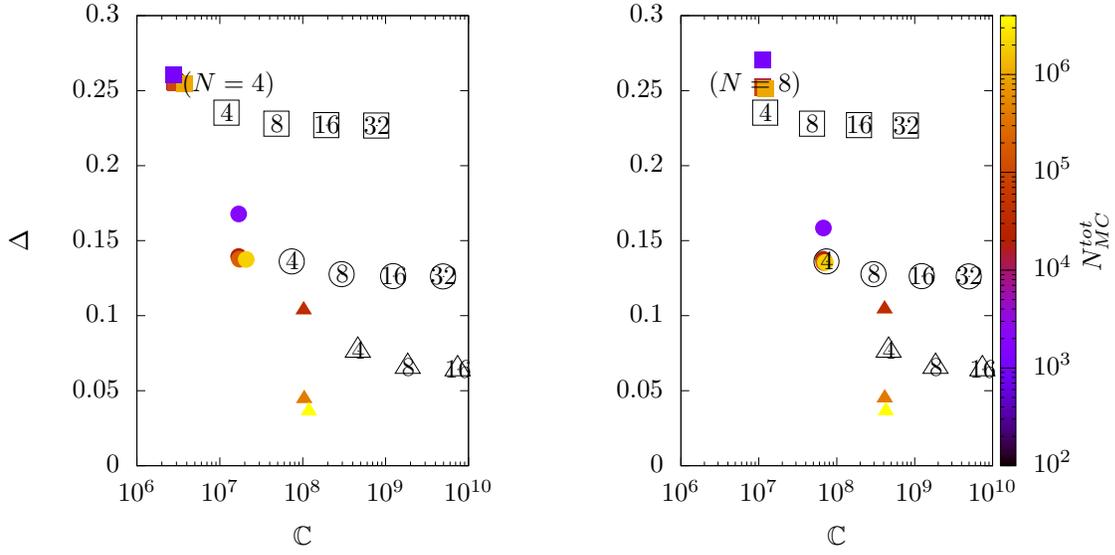

Figure 7: The relative $L^2$-difference $\Delta$ vs. complexity $\mathbb{C}$ for the scalar flux $\Phi$ in the lattice problem, using the hybrid method (filled, colored markers) and the monolithic $S_N$ method (empty, green markers). Points that are down and to the left are more efficient. The hybrid-solver was run for $N = 4$ (left) and $N = 8$ (right). The shape of data-points indicates different spatial resolution (square: $N_{\boldsymbol{x}} = 56$, circle: $N_{\boldsymbol{x}} = 112$, triangle: $N_{\boldsymbol{x}} = 224$). Coloring of the Hybrid-data points corresponds to the total number of particles per time step $N_{MC}^{tot}$ according to the colorbar. The $S_N$ method was run for $N = 4, 8, 16, 32$. Each DG-data point is assigned a numerical label according to its value of $N$. The formula for $\Delta$ is given in (48) which the complexity $\mathbb{C}$ is given by (49a) for the monolithic $S_N$ method and by (49b) for the hybrid.

*3.3. The linearized hohlraum problem*

In the linearized hohlraum problem [9], nonlinear coupling between particles and the material medium is approximated in a linear way by adjusting the absorption and scattering cross-sections according to the expected material temperature profile of the nonlinear problem [8]. The geometry of the setup along with the material parameters can be found in Figure 8. The domain is $X = [0, 1.3] \times [0, 1.3]$, and the initial condition is identically zero everywhere. For boundary conditions, we assume a constant influx from the left side of the domain, i.e. $\Psi(x = 0, y, \mathbf{\Omega}, t) = 1$ for $\Omega_x > 0$. As discussed in the appendix, this boundary condition can be treated as a surface source, modeled by setting $\Omega_x = \sqrt{\xi}$, where $\xi \sim U([0, 1])$ is sampled uniformly on $[0, 1]$. The spatial distribution along the boundary is sampled uniformly.

We again perform $S_N$ and hybrid runs with varying spatial and angular resolution. The spatial domain is subdivided into equal $N_x \times N_x$ square cells with $N_x \in \{52, 104, 208\}$. For the $S_N$-runs we use $N \in \{4, 8, 16, 32\}$, resulting in $N_\Omega = N^2$ ordinates on the northern hemisphere of $\mathbb{S}^2$. The collided part of the hybrid algorithm employs an $S_N$ method with $N \in \{4, 8\}$. In the hybrid method the number of particles is also changed between runs, but the killing weight remains fixed at $w_{\text{kill}} = 10^{-15}$. The number of particles newly inserted into the system every time step is $2 \times 10^k$ for $k \in \{2, 3, 4, 5, 6\}$. All runs are performed to a final time of $t_{\text{final}} = 2.6$ and the CFL is kept fixed at 52. The reference solution is a $S_{96}$-$S_{16}$ hybrid using a $T_N$ quadrature in angle, a third-order DG method in space on a $448 \times 448$ grid, and a defect correction time integrator [10].

In Figure 9, we show densities of a few select runs, calculated using $S_N$ and hybrid methods. The log of the corresponding relative $L_2$-errors are depicted in Figure 10. As before, the $S_N$ solutions suffer from ray-effects that are marginally reduced as the number of angle increases. Meanwhile, most of the disparities between the reference and hybrid solutions can be attributed to stochastic noise. The hybrid has a mix of ray effects from the collided equation solve and particle tracks from the uncollided equation solve on the backside of the hohlraum However, these errors here are on the order of $10^{-2}$-$10^{-3}$, which is much smaller than the errors in othe back of the domain.

Detailed comparisons between the relative error against the computational complexity are depicted in Figure 11. As before, increasing the angular resolution in the collided equation does not benefit the accuracy of the hybrid method. Increasing particles also has less effect than in the previous problems. The $S_N$ solutions benefit most from finer spatial resolution, while the angular resolution does not matter as much. Spatial resolution also plays the biggest role for the hybrid. We do observe that for small particle counts (the purple points in the figure) increasing resolution can actually increase the error. This is explained by the fact that an under-sampled MC calculation does not benefit from more spatial resolution. Nevertheless, for a fixed spatial resolution, we do observe an improvement in the error when $N_p$ is increased.

Overall the hybrid runs produce solutions with comparable or better accuracy than monolithic $S_N$ runs with the same spatial resolution. Hybrid runs using $S_8$ for the collided equation achieve improved accuracy at nearly the same complexity while runs using $S_4$ for the collided equation achieved improved accuracy with even less complexity.

## 4. Conclusion and Discussion

In this work, we have presented a collision-based hybrid method that uses a Monte Carlo method for the uncollided solution and a discrete ordinate discretization for the collided solution. This combination of methods was originally proposed for the first collision source strategy used in [3] in the steady-state setting. Thus this work can be considered as an extension to the time-dependent setting that requires a remap procedure after every time step.

Experimental simulations have been performed on three standard benchmarks. For each benchmark, the results demonstrate that the hybrid method is more efficient, in the sense that it achieves greater accuracy with the comparable or less complexity or is less complexity with comparable or greater accuracy. Here complexity is a measure of how many unknowns are updated during particle moves for Monte Carlo or sweeping iterations for discrete ordinates.

This work has concentrated on single-energy particle transport problems. However recent work has shown that when considering energy-dependent problems, more opportunities for hybridization arise when considering fully deterministic hybrids [40]. In those results it was shown that low-resolution in energy can be used for the collided solution as well as low-resolution in angle. With the introduction of Monte Carlo, new opportunities arise. For example, continuous energy cross-sections could be used in the uncollided portion.

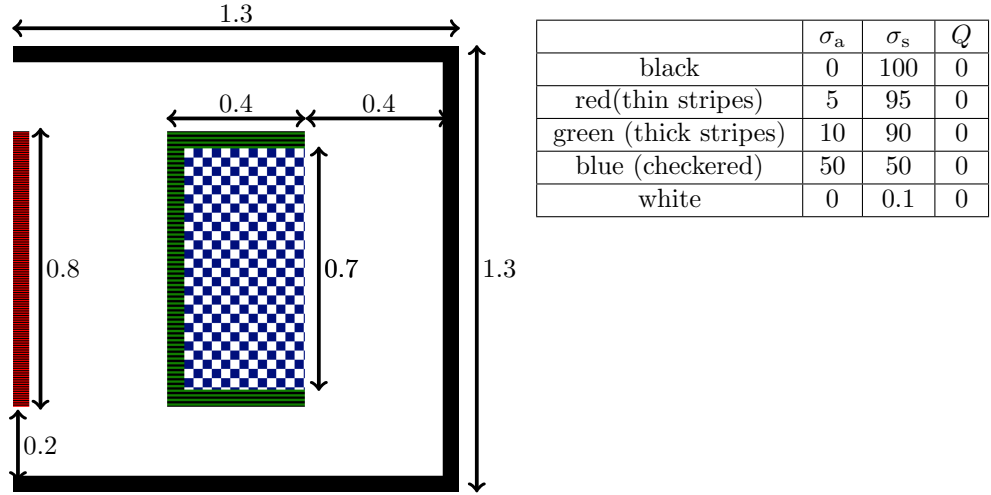

Figure 8: Geometric layout and table of material parameters for the Hohlraum Problem. All walls have a thickness of 0.05.

This could be important to treat resonances in neutron transport problems, but investigation is needed to quantify any benefits from this approach. In addition future investigations should be made regarding the use of Monte Carlo techniques inside the high-order time accuracy methods developed for hybrid problems in [9, 10].

## 5. Acknowledgments

J.K. gratefully acknowledges support from the 2022 National Science Foundation Mathematical Sciences Graduate Internship to conduct this research at Oak Ridge National Laboratory.

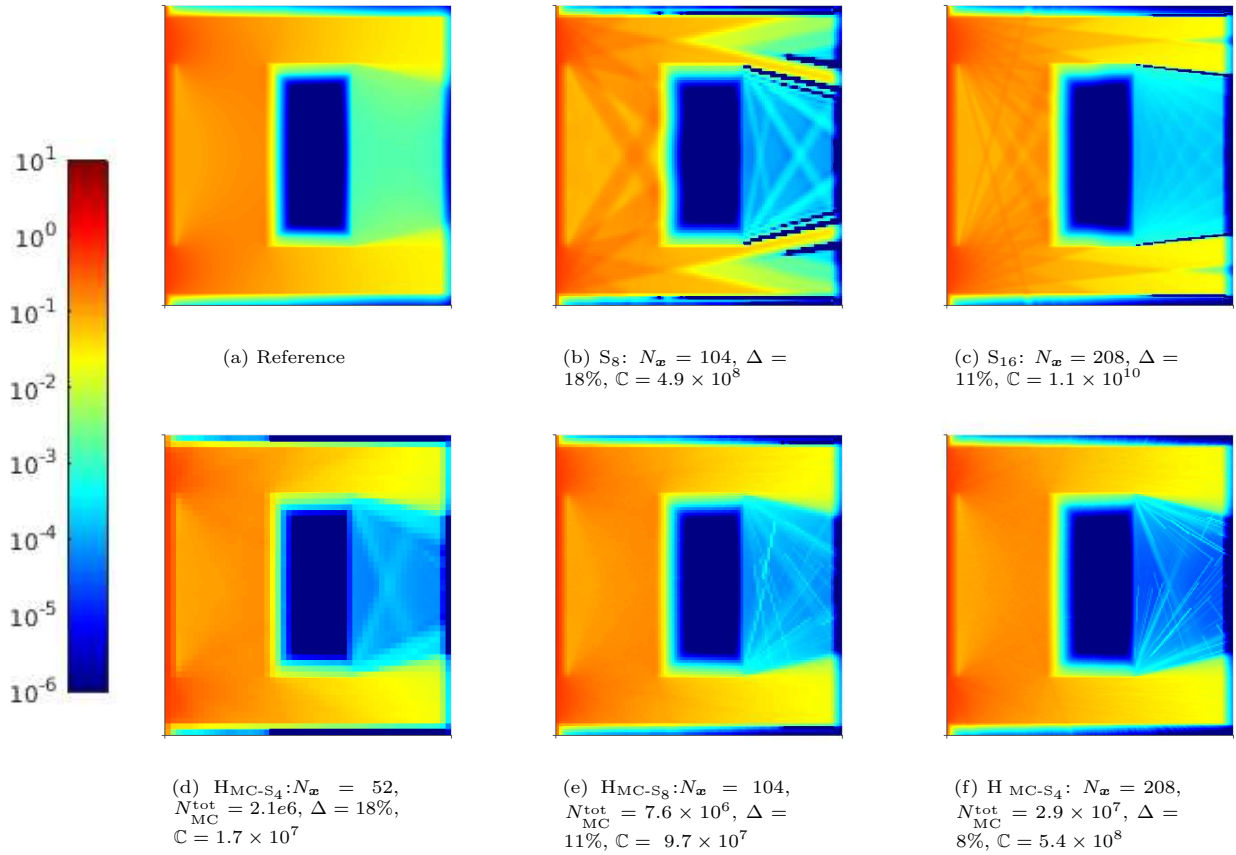

Figure 9: Numerical solutions to the hohlraum problem at $t = 2.6$ with CFL 52. Each numerical solution is characterized by a relative $L^2$ difference $\Delta$ with respect to the reference, defined in (48), and a complexity $\mathbb{C}$, defined in (49a) for the monolithic $S_N$ method and (49b) for the hybrid.

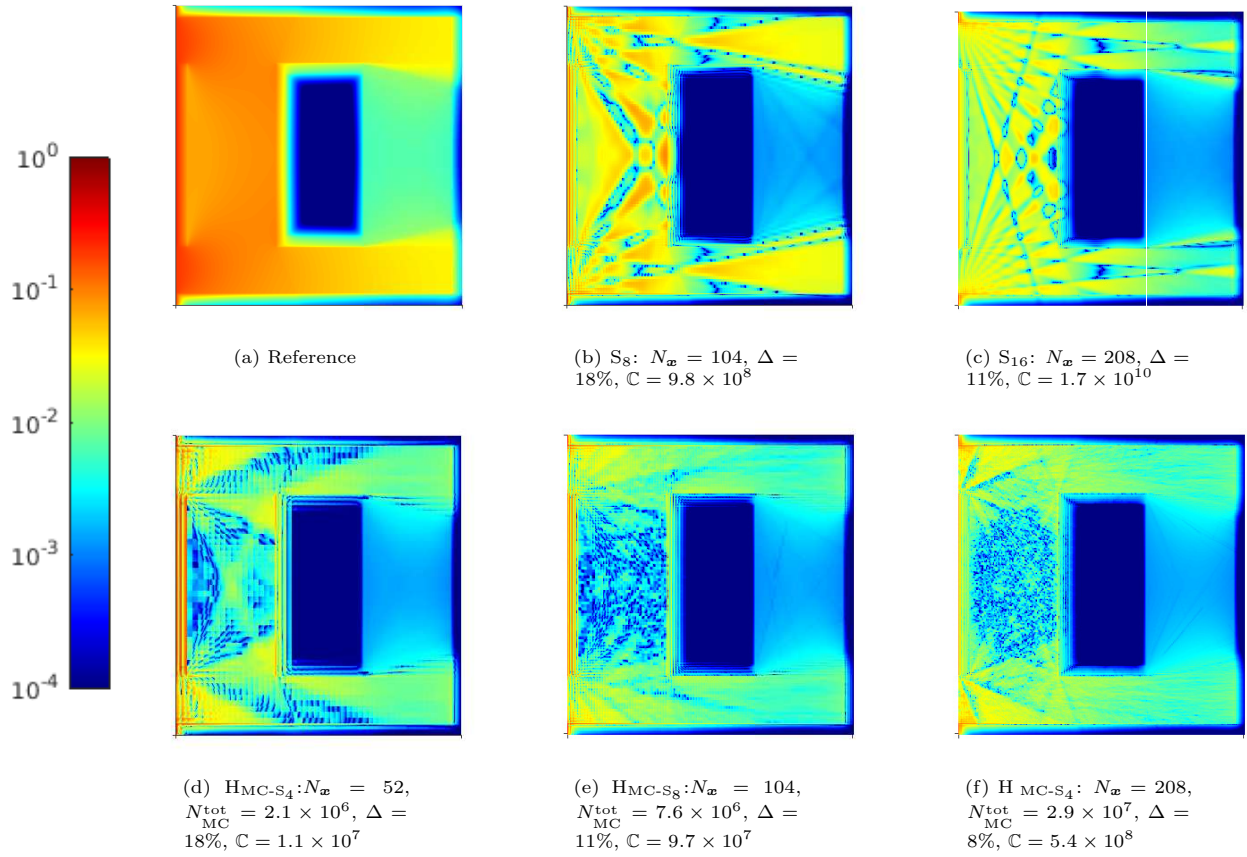

Figure 10: Hohlraum problem: Absolute difference between analytical solution to the line source problem and various numerical solutions at $t = 2.6$ with CFL 52. Each numerical solution is characterized by a relative $L^2$ difference $\Delta$ with respect to the reference, defined in (48), and a complexity $\mathbb{C}$, defined in (49a) for the monolithic $S_N$ method and (49b) for the hybrid.

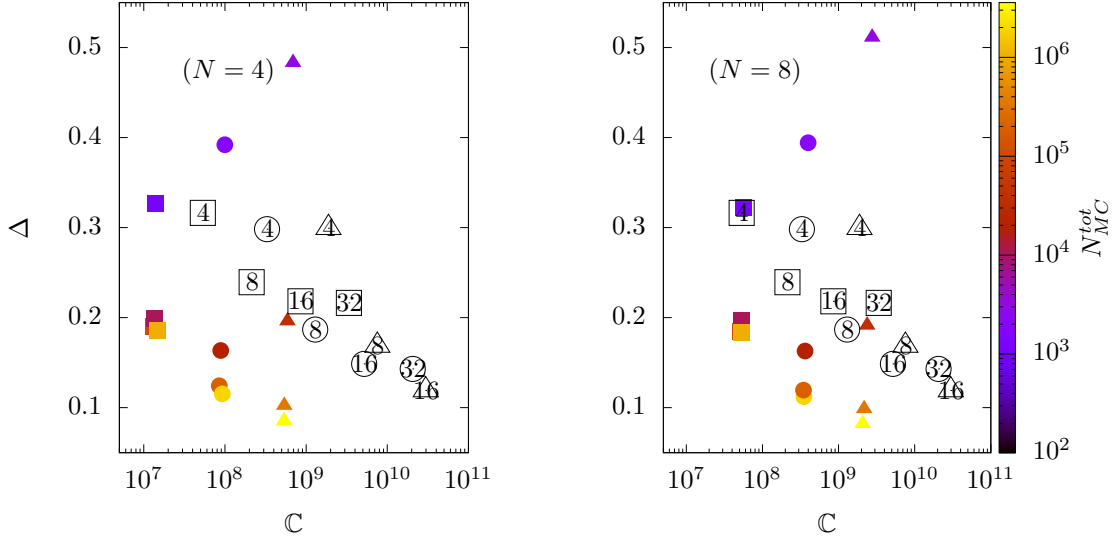

Figure 11: The relative $L^2$-difference $\Delta$ vs. complexity $\mathbb{C}$ for the scalar flux $\Phi$ in the hohlraum problem, using the hybrid method (filled, colored markers) and the monolithic $S_N$ method (empty, green markers). Points that are down and to the left are more efficient. The hybrid-solver was run for $N = 4$ (left) and $N = 8$ (right). The shape of data-points indicates different spatial resolution (square: $N_{\boldsymbol{x}} = 52$, circle: $N_{\boldsymbol{x}} = 104$, triangle: $N_{\boldsymbol{x}} = 208$). Coloring of the Hybrid-data points corresponds to the total number of particles per time step $N_{MC}^{tot}$ according to the colorbar. The $S_N$ method was run for $N = 4, 8, 16, 32$. Each DG-data point is assigned a numerical label according to its value of $N$. The formula for $\Delta$ is given in (48) which the complexity $\mathbb{C}$ is given by (49a) for the monolithic $S_N$ method and by (49b) for the hybrid.

### Appendix A. Boundary conditions for the hohlraum problem

Unlike the line source and lattice problems, the linearized hohlraum problem involves non-zero boundary conditions. A Monte Carlo implementation of this boundary condition is stated in [15]; here we present a derivation of the approach that is used.

In the hohlraum problem, it is assumed that $\Psi(x, y, t, \Omega) = 1$ for $x = 0$ and $\Omega_x > 0$, i.e., a constant flux of 1 along the left boundary is assumed for each incoming direction. To model this with Monte Carlo, we assume that this flux is due to a source $s_b$ (see (37c)) located on an infinitesimal slab just left of the boundary.

Consider first a finite slab $S_a = \{(x, y) \in [-a, 0] \times [0, 1.3]\}$, where $a > 0$. We assume that $\sigma_{\mathrm{a}} = \sigma_{\mathrm{s}} = 0$ on $S$, that the source $s_b(x, y, \boldsymbol{\Omega}; a) = s_b(x, \boldsymbol{\Omega}; a)$ is independent of $y$ and $t$, and that

$$\hat{s}_b(\boldsymbol{\Omega}) := \int_{-a}^{0} s_b(x, \boldsymbol{\Omega}; a) dx \qquad (A.1)$$

is independent of $a$. Thus, $s_b(\mathbf{x}, \boldsymbol{\Omega}, t) \to \hat{s}_b(\boldsymbol{\Omega})\delta(x)$ as $a \to 0$.

To determine $\hat{s}_b$, we assume that $\Psi$ is independent of $y$ and $t$ on $S$ and satisfies the steady-state equation

$$\Omega_x \Psi_x(x, y, \boldsymbol{\Omega}) = s_b, \qquad (x, y) \in S, \quad \boldsymbol{\Omega} \in \mathbb{S}^2, \qquad (A.2a)$$

$$\Psi(-a, y, \boldsymbol{\Omega}) = 0, \qquad y \in [0, 1.3], \quad \Omega_x > 0. \qquad (A.2b)$$

This formulation is consistent with (34), given the assumptions made on $\Psi$, $s_b$, and the material cross-sections. Integrating (A.2b) with respect to $x$ and applying the boundary condition in (A.2b) gives

$$\hat{s}_b(\boldsymbol{\Omega}) = \Omega_x, \qquad \Omega_x > 0. \qquad (A.3)$$

Hence $s_b = \Omega_x \delta(x)$.

Finally, to sample particles according to the probability density $\rho(\Omega_x) \propto \hat{s}_b$ for $\Omega_x > 0$, we compute the cumulative distribution function (CDF):

$$F(\Omega_x) = \int_0^{\Omega_x} 2\mu d\mu = \Omega_x^2. \tag{A.4}$$

According to the fundamental theorem of simulation [30, pp. 19-22], the correct angular distribution can be sampled by generating uniform variables $u \in [0,1]$ and setting $\Omega_x = F^{-1}(u) = \sqrt{u}$.